\documentclass[11pt,twoside,leqno]{aomamlt2e} 
  \pageno{777}
\received{February 5, 2003}

\newtheorem{thm}{Theorem}[section]
\newtheorem{lem}[thm]{Lemma}
\newtheorem{prop}[thm]{Proposition}
\newtheorem{cor}[thm]{Corollary}

\theorembodyfont{\upshape}
\newtheorem{defn}[thm]{{\it Definition}}
\newtheorem{sct}[thm]{ }
\newtheorem{rem}[thm]{{\it Remark}}
\numberwithin{equation}{section}
\numberwithin{figure}{section}

\def\longr{{{\longrightarrow}}}
\def\defo{{{\mathrm{Def}}}}

\def\Kern{{{\mathrm{Ker}}}}
\def\Im{{{\mathrm{Im}}}}
\def\Ext{{{\mathrm{Ext}}}}
\def\Tor{{{\mathrm{Tor}}}}
\def\Der{{{\mathrm{Der}}}}
\def\Hom{{{\mathrm{Hom}}}}
\def\Spec{{{\mathrm{Spec}}}}

\def\map{{{\mathrm{map}}}}
\def\mm{{{\mathfrak{m}}}}
\def\AA{{\mathbb A}}

\def\FF{{\mathbb F}}
\def\GG{{\mathbb G}}
\def\WW{{\mathbb W}}
\def\SS{{\mathbb S}}
\def\ZZ{{\mathbb Z}}

\def\cF{{\cal F}}
\def\cG{{\cal G}}
\def\cK{{\cal K}}
\def\cL{{\cal L}}
\def\cO{{\cal O}}
\def\Chi{{\chi}}
\def\FGL{{\cF\cG\cL}}

\def\RINGS{{{\rm Rings}_c}}
\def\colim{{{\mathrm{colim}}}}
\def\hocolim{{{\mathrm{hocolim}}}}

\def\holim{{\mathrm{holim}}}

\def\defeq{\overset{\mathrm{def}}=}
\def\Gal{{{\mathrm{Gal}}}}

\renewcommand{\theequation}{\arabic{section}.\arabic{equation}}

\begin{document}
\currannalsline{162}{2005} 

\title{A resolution of the $K(2)$-local sphere\\
 at the prime $3$}

 \acknowledgements{The first
author and fourth authors were partially supported by the National Science
Foundation (USA). The authors would like to thank (in alphabetical order) 
MPI at Bonn, Northwestern University, the Research in Pairs 
Program at Oberwolfach, the University of Heidelberg and 
Universit\'e Louis Pasteur at Strasbourg, 
for providing them with the opportunity to work together.}
 
\twoauthors{P.\ Goerss, H.-W.  Henn, M.\ Mahowald,}{C.  Rezk}
\shortname{P. Goerss, H.-W. Henn, M. Mahowald, and C. Rezk}
\shorttitle{A resolution of the $K(2)$-local sphere at the prime $3$}
 \institution{Northwestern University, Evanston, IL
\\
\email{pgoerss@math.northwestern.edu}\\
\vglue-9pt
Institut de Recherche Math\'ematique Avanc\'ee,
Universit\'e Louis Pasteur,\\ Strasbourg, France
\\
\email{henn@math.u-strasbg.fr}
\\
\vglue-9pt
Northwestern University, Evanston, IL
\\
\email{mark@math.northwestern.edu}
\\
\vglue-9pt
University of Illinois at Urbana-Champaign, Urbana, IL
\\
\email{rezk@math.uiuc.edu}}

 \vglue12pt

\centerline{\bf Abstract} 
\vglue12pt
We develop a framework for displaying
the stable homotopy theory of the sphere, at least after localization 
at the second Morava $K$-theory $K(2)$. At the
prime $3$, we write the spectrum $L_{K(2)}S^0$ as the inverse
limit of a tower of fibrations with four layers. The successive fibers are
of the form $E_2^{hF}$ where $F$ is a finite subgroup of the Morava
stabilizer group and $E_2$ is the second Morava or Lubin-Tate
homology theory. We give explicit calculation of the homotopy groups
of these fibers. The case $n=2$ at $p=3$ represents the edge
of our current knowledge: $n=1$ is classical and at $n=2$, the prime $3$
is the largest prime where the Morava stabilizer group has a $p$-torsion
subgroup, so that the homotopy theory is not entirely algebraic.

\vglue16pt

The problem of understanding the homotopy groups of spheres has been
central to algebraic topology ever since the field emerged as a 
distinct area of mathematics. A period of calculation beginning with Serre's
computation of the cohomology of Eilenberg-MacLane spaces and the
advent of the Adams spectral sequence culminated, in the late 1970s,
with the work of Miller, Ravenel, and Wilson on periodic phenomena
in the homotopy groups of spheres and Ravenel's nilpotence
conjectures. The solutions to most of these conjectures by
Devinatz, Hopkins, and Smith in the middle 1980s established the primacy
of the ``chromatic'' point of view and there followed a period in which
the community absorbed these results and extended the qualitative picture
of stable homotopy theory. Computations passed from center stage, to
some extent, although there has been steady work in the wings --
most notably by Shimomura and his coworkers, and Ravenel, 
and more lately by Hopkins and his coauthors in their work on
topological modular forms. The amount of interest generated by
this last work suggests that we may be entering a period of
renewed focus on computations.

In a nutshell, the chromatic point of view is based on the observation
that much of the structure of stable homotopy theory is controlled by the
algebraic geometry of formal groups. The underlying geometric object
is the moduli stack of formal groups. Much of what can be
proved and conjectured about stable homotopy theory arises
from the study of this stack, its stratifications, and the theory
of its quasi-coherent sheaves. See for example, the table in
Section 2 of \cite{HG}. 

The output we need from this geometry consists of two distinct pieces
of data. First, the chromatic convergence theorem of \cite[\S 8.6]{Nil} says the following. Fix a prime $p$
and let $E(n)_\ast$,
$n \geq 0$ be the Johnson-Wilson homology theories and let $L_n$ be
localization with respect to $E(n)_\ast$. Then there are natural
maps $L_nX \to L_{n-1}X$ for all spectra $X$, and if $X$
is a $p$-local finite spectrum, then the natural map
$$
X \longr\ \holim L_n X
$$
is a weak equivalence. 

Second, the maps $L_n X \to L_{n-1}X$ fit into a good
fiber square. Let $K(n)_\ast$ denote the $n$-th Morava $K$-theory.
Then there is a natural commutative diagram 
\begin{equation}\label{fiber-square}
\xymatrix{
L_n X \rto \dto & L_{K(n)}X \dto\\
L_{n-1}X \rto & L_{n-1}L_{K(n)}X\\
}
\end{equation}
which for any spectrum $X$ is a homotopy pull-back square. 
It is somewhat difficult to find this result in
the literature; it is implicit in \cite{chromatic}.

Thus, if $X$ is a  $p$-local finite spectrum, the basic building blocks for
the homotopy type of $X$ are the Morava $K$-theory localizations $L_{K(n)}X$.

Both the chromatic convergence theorem and the fiber square of
(\ref{fiber-square}) can be viewed as analogues of phenomena familiar
in algebraic geometry. For example, the fibre square can be thought of as 
an analogue of a Mayer-Vietoris situation for a formal neighborhood of a 
closed subscheme and its open complement (see \cite{BL}).  
The chromatic convergence theorem can be thought of as a result
which determines  
what happens on a variety $S$ with a nested sequence of 
closed sub-schemes $S_n$ of codimension $n$ by what happens on the open 
subvarieties $U_n=S-S_n$ (See \cite[\S IV.3]{Hart}, for example.) 
This analogy can be made precise using the moduli stack of $p$-typical
formal group laws for $S$ and, for $S_n$, the substack which classifies
formal groups of height at least $n$. Again see \cite{HG}; also, see
\cite{Pribble} for more details. 
 
In this paper, we will write (for $p=3$) the $K(2)$-local stable sphere as
a very small 
homotopy inverse limit of spectra with computable and computed
homotopy groups. Specifying a Morava $K$-theory always means fixing
a prime $p$ and a formal group law of height $n$; we unapologetically focus
on the case $p=3$ and $n=2$ because this is at the edge
of our current knowledge. The homotopy type and homotopy groups
for $L_{K(1)}S^0$ are well understood at all primes and are
intimately connected with the $J$-homomorphism; indeed, this
calculation was one of the highlights of the computational period
of the 1960s. If $n = 2$ and $p > 3$, the Adams-Novikov spectral
sequence (of which more is said below) calculating $\pi_\ast L_{K(2)}S^0$ 
collapses and cannot have extensions; hence, the problem becomes
algebraic, although not easy. Compare \cite{ShimY}.

It should be noticed immediately that for $n=2$ and $p=3$ 
there has been a great deal of calculations of the homotopy groups of 
$L_{K(2)}S^0$ and closely related spectra, 
most notably by Shimomura and his coauthors. 
(See, for example, \cite{Shimv}, \cite{Shim} and \cite{ShimW}.) 
One aim of this paper is to provide a conceptual framework for 
organizing those results and produce further advances.

The $K(n)$-local category of spectra is governed by a
homology theory built from the Lubin-Tate (or Morava) theory $E_n$. This
is a commutative ring spectrum with coefficient ring
$$
(E_n)_\ast = W(\FF_{p^n})[[u_1,\dots ,u_{n-1}]][u^{\pm 1}]
$$
with the power series ring over the Witt vectors in degree $0$ and
the degree of $u$ equal to $-2$. The ring
$$
(E_n)_0 = W(\FF_{p^n})[[u_1,\dots ,u_{n-1}]]
$$
is a complete local ring with residue field $\FF_{p^n}$. It is
one of the rings constructed by Lubin and Tate in their study
of deformations for formal group laws over fields of
characteristic $p$. See \cite{LT}.

As the notation indicates, $E_n$ is closely related to the
Johnson-Wilson spectrum $E(n)$ mentioned above.

The homology theory $(E_n)_\ast$ is a complex-oriented theory
and the formal group law over $(E_n)_\ast$ is a universal deformation
of the Honda formal group law $\Gamma_n$ of height $n$ over the field
$\FF_{p^n}$ with $p^n$ elements. 
(Other choices of formal group laws of height $n$  are possible, but all
yield essentially the same results. The choice of $\Gamma_n$ is only made 
to be explicit; it is the usual formal group law associated 
by homotopy theorists to Morava $K$-theory.) 
Lubin-Tate theory implies that the graded ring $(E_n)_\ast$ supports
an action by the group
$$
\GG_n = \hbox{Aut}(\Gamma_n) \rtimes {\rm Gal}(\FF_{p^n}/\FF_p).
$$
The group $\hbox{Aut}(\Gamma_n)$ of automorphisms of the formal group
law $\Gamma_n$ is also known as the Morava stabilizer group and
will be denoted by $\SS_n$. The Hopkins-Miller theorem (see \cite{Rezk})
says, among other things, that we can lift this action to an
action on the spectrum $E_n$ itself.
There is an Adams-Novikov spectral
sequence
$$
E_2^{s,t}:=H^s(\SS_n,(E_n)_t)^{{\rm Gal}(\FF_{p^n}/\FF_p)} 
\Longrightarrow \pi_{t-s}L_{K(n)}S^0.
$$
(See \cite{HMS} for a basic description.) The group $\GG_n$ is a profinite
group and it acts continuously on $(E_n)_\ast$. The cohomology here is
continuous cohomology. We note that by \cite{DH1} $L_{K(n)}S^0$ can be 
identified with the homotopy fixed point spectrum $E_n^{h\GG_n}$ 
and the Adams-Novikov spectral sequence can 
be interpreted as a homotopy fixed point spectral sequence. 

The qualitative behaviour of this 
spectral sequence depends very much on qualitative 
cohomological properties of the group 
$\SS_n$, in particular on its cohomological dimension. 
This in turn depends very much on $n$ and $p$.  

If $p-1$ does not divide $n$ (for example,
if $n < p-1$) then the\break $p$-Sylow subgroup of $\SS_n$ is 
of cohomological dimension $n^2$. 
Furthermore, if\break $n^2 < 2 p-1$ (for example, if $n=2$ and $p >3$) 
then this spectral sequence is sparse enough so that there
can be no differentials or extensions.

However, if $p-1$ divides $n$, then 
the cohomological dimension of $\SS_n$ is infinite and 
the Adams-Novikov spectral sequence has a more complicated behaviour. 
The reason for infinite cohomological dimension is the existence of 
elements of order $p$ in $\SS_n$. 
However, in this case at least the virtual cohomological dimension remains 
finite, in other words there are finite index subgroups with finite 
cohomological dimension. 
In terms of resolutions of the trivial module $\ZZ_p$, this means 
that while there are no projective resolutions of the trivial 
$\SS_n$-module $\ZZ_p$ of finite length, one might still hope that there 
exist ``resolutions'' of $\ZZ_p$ 
of finite length in which the individual modules are direct 
sums of modules which are permutation modules 
of the form $\ZZ_p[[\GG_2/F]]$ where $F$ is a finite subgroup of $\GG_n$. 
Note that in the case of a discrete group which acts properly and 
cellularly on a finite dimensional contractible space $X$ such a 
``resolution'' is provided by the complex of cellular chains on $X$.

This phenomenon is already visible for $n=1$ in which case 
$\GG_1=\SS_1$ can be identified with $\ZZ_p^{\times}$, 
the units in the $p$-adic integers. Thus $\GG_1\cong \ZZ_p\times C_{p-1}$ 
if $p$ is odd while $\GG_1\cong \ZZ_2\times C_2$ if $p=2$. 
In both cases there is a short exact sequence 
$$
0\to \ZZ_p[[\GG_1/F]]\to \ZZ_p[[\GG_1/F]]\to \ZZ_p\to 0
$$
of continuous $\GG_1$-modules (where $F$ is 
the maximal finite subgroup of $\GG_1$).   
If $p$ is odd this sequence is a projective 
resolution of the trivial module while for $p=2$ it is only a resolution 
by permutation modules. These resolutions are the algebraic 
analogues of the fibrations (see \cite{HMS}) 
\begin{equation}\label{LK1}
L_{K(1)}S^0\simeq E_1^{h\GG_1}\to E_1^{hF}\to E_1^{hF}\ .
\end{equation}
We note that $p$-adic complex $K$-theory 
$K\ZZ_p$ is in fact a model for $E_1$,  
the homotopy fixed points $E_1^{hC_{2}}$ can be identified 
with $2$-adic real $K$-theory $KO\ZZ_2$ if $p=2$ 
and $E_1^{hC_{p-1}}$ is the 
Adams summand of $K\ZZ_p$ if $p$ is odd,  so that the fibration of 
(\ref{LK1}) indeed agrees with that of \cite{HMS}.

In this paper we produce a resolution of the trivial module $\ZZ_p$ by 
(direct summands of) permutation modules in the case $n=2$ and $p=3$  
and we use it to build $L_{K(2)}S^0$ as the top of a finite tower of fibrations
where the fibers are (suspensions of) spectra of the
form $E_2^{hF}$ where $F\subseteq \GG_2$ is a finite subgroup.

In fact, if $n=2$ and $p=3$, only two subgroups
appear. The first is a subgroup $G_{24} \subseteq \GG_2$; this is
a finite subgroup of order $24$ containing a normal cyclic subgroup $C_3$ with 
quotient $G_{24}/C_3$ isomorphic to the quaternion group $Q_8$ of order $8$.
The other group is the semidihedral group $SD_{16}$ of order $16$. 
The two spectra we will see, then, are $E_2^{hG_{24}}$ and 
$E_2^{hSD_{16}}$. 

The discussion of these and related subgroups of $\GG_2$ occurs in Section
1 (see \ref{rem0.1} and  \ref{rem0.2}). 
The homotopy groups of these spectra are known. 
We will review the calculation in Section 3. 

Our main result can be stated as follows (see Theorems 
\ref{thm5.4} and  \ref{thm5.5}).

\begin{thm}
There is a sequence of maps between spectra 
\begin{eqnarray*}
L_{K(2)}S^0 \to E_2^{hG_{24}} &\to &\Sigma^8E_2^{hSD_{16}} 
\vee E_2^{hG_{24}} 
 \to 
\Sigma^8E_2^{hSD_{16}} \vee \Sigma^{40}E_2^{hSD_{16}}  
\\
&\to& \Sigma^{40}E_2^{hSD_{16}} 
\vee \Sigma^{48}E_2^{hG_{24}} \to \Sigma^{48}E_2^{hG_{24}}
\end{eqnarray*}
with the property that the composite of any two successive maps 
is zero and all
possible Toda brackets are zero modulo indeterminacy. 
\end{thm}

Because the Toda brackets vanish, this ``resolution'' can be refined
to a tower of spectra with $L_{K(2)}S^0$ at the top. The
precise result is given in Theorem~\ref{thm5.6}. There are many
curious features of this resolution, of which we note here
only two. First, this is not an Adams resolution for $E_2$, as
the spectra $E_2^{hF}$ are not $E_2$-injective, at least if
$3$ divides the order of $F$. Second, there is a certain superficial duality 
to the resolution which should somehow be explained by the fact that 
$\SS_n$ is a virtual Poincar\'e
duality group, but we do not know how to make this
thought precise. 

As mentioned above, this result can be used to organize the already 
existing and very complicated calculations of Shimomura (\cite{Shim}, 
\cite {ShimW}) and it also suggests an independent approach to 
these calculations. Other applications would be to the study of
Hopkins's Picard group (see \cite{HMS}) of $K(2)$-local invertible spectra. 

Our method is by brute force. The hard work is really in Section 4, where
we use the calculations of \cite{Henn} in an essential way to
produce the short resolution of the trivial $\GG_2$-module $\ZZ_3$
by (summands of) permutation modules of the form $\ZZ_3[[\GG_2/F]]$
where $F$ is finite (see Theorem 
\ref{thm4.1} and Corollary~\ref{cor4.2}). 
In Section 2, we calculate the homotopy type of
the function spectra $F(E^{hH_1},E^{hH_2})$ 
if $H_1$ is a closed and 
$H_2$ a finite subgroup of $\GG_n$; this will allow
us to construct 
\newpage
\noindent
the required maps between these spectra and to 
make the Toda bracket calcula-
\-
tions. Here the work of
\cite{DH1} is crucial. These calculations also explain the
role of 
\linebreak
the suspension by 48 which is really 
a homotopy theoretic phenomenon while 
\linebreak
the 
other suspensions can be explained in terms of the algebraic resolution 
\linebreak
constructed in Section 4.

\section{Lubin-Tate theory and the Morava stabilizer group}

The purpose of this section is to give a summary of what
we will need about deformations of formal group laws over perfect fields.
The primary point of this section is to establish notation
and to run through some of the standard algebra needed to come
to terms with the $K(n)$-local stable homotopy category.

Fix a perfect field $k$ of characteristic $p$ and a formal group law
$\Gamma$ over $k$. A {\it deformation} of $\Gamma$ to a 
complete local ring $A$ (with maximal ideal $\mm$) is a pair $(G,i)$ 
where $G$ is a formal group law over $A$, 
$i:k \to A/\mm$ is a morphism of fields
and one requires $i_\ast\Gamma = \pi_\ast G$, where $\pi:A \to A/\mm$
is the quotient map. Two such deformations $(G,i)$ and $(H,j)$ are
$\star$-isomorphic if there is an isomorphism $f:G \to H$ of
formal group laws which reduces to the identity modulo $\mm$. 
Write $\defo_\Gamma(A)$ for the set of $\star$-isomorphism
classes of deformations of $\Gamma$ over $A$.

A common abuse of notation is to write $G$ for the deformation
$(G,i)$; $i$ is to be understood from the context.

Now suppose the height of $\Gamma$ is finite.
Then the theorem of Lubin and Tate \cite{LT}
says that the functor $A \mapsto \defo_\Gamma(A)$ is representable.
Indeed let
\begin{equation}\label{eq0.1}
E(\Gamma,k) = W(k)[[u_1,\dots ,u_{n-1}]]
\end{equation}
where $W(k)$ denotes the Witt vectors on $k$ and $n$ is the height of 
$\Gamma$. This is a complete
local ring with maximal ideal $\mm = (p,u_1,\dots ,u_{n-1})$ and there
is a canonical isomorphism $q: k \cong E(\Gamma,k)/\mm$. Then Lubin and Tate
prove there is a deformation $(G,q)$ of $\Gamma$ over $E(\Gamma,k)$ 
so that the natural map
\begin{equation}\label{eq0.2}
\Hom_c (E(\Gamma,k),A) \to \defo_\Gamma(A)
\end{equation}
sending a continuous map $f:E(\Gamma,k) \to A$ to $(f_\ast G, \bar{f}q)$
(where $\bar{f}$ is the map on residue fields induced by $f$) 
is an isomorphism. Continuous maps here are very simple: they are the
{\it local} maps; that is, we need only
require that $f(\mm)$ be contained in the maximal ideal of $A$.
Furthermore, if two deformations are $\star$-isomorphic, then the 
$\star$-isomorphism between them is unique.

We would  like to now turn the assignment $(\Gamma,k) \mapsto E(\Gamma,k)$
into a functor. For this we introduce the  category $\FGL_n$ of height
$n$ formal group laws over perfect fields. The objects are
pairs $(\Gamma,k)$ where $\Gamma$ is of height $n$. A morphism
$$
(f,j):(\Gamma_1,k_1) \to (\Gamma_2,k_2)
$$
is a homomorphism of fields $j:k_1 \to k_2$ and an isomorphism of
formal group laws $f:j^\ast\Gamma_1 \to \Gamma_2$.

Let $(f,j)$ be such a morphism and let 
$G_1$ and $G_2$ be the fixed universal deformations over $E(\Gamma_1,k_1)$
and $E(\Gamma_2,k_2)$ respectively. 
If $\widetilde{f} \in E(\Gamma_2,k_2)[[x]]$ is any lift
of $f \in k_2[[x]]$, then we can define a formal group
law $H$ over $ E(\Gamma_2,k_2)$\break\vskip-12pt\noindent  by requiring that $\widetilde{f}:H \to G_2$
is an isomorphism. Then the pair $(H,j)$ is
a deformation of $\Gamma_1$, hence we get a homomorphism 
$E(\Gamma_1,k_1) \to E(\Gamma_2,k_2)$ 
classifying the $\star$-isomorphism class of $H$ --
which, one easily checks, is independent of the lift $\widetilde{f}$.
Thus if $\RINGS$ is the category of complete local rings and local
homomorphisms, we get a functor
$$
E(\cdot,\cdot):\FGL_n \longrightarrow \RINGS.
$$
In particular, note that any morphism in $\FGL_n$ from a pair $(\Gamma,k)$
to itself is an isomorphism. The automorphism group of 
$(\Gamma,k)$ in $\FGL_n$ is the ``big'' Morava stabilizer group
of the formal group law; it contains the subgroup of elements
$(f,{\rm id}_k)$. This formal group law and hence also its automorphism group 
is determined up to isomorphism 
by the height of $\Gamma$ if $k$ is separably closed. 

Specifically, let $\Gamma$ be the Honda formal group law over
$\FF_{p^n}$; thus the $p$-series of $\Gamma$ is 
$$
[p](x) = x^{p^n}.
$$
From this formula it immediately follows that any automorphism
$f:\Gamma \to \Gamma$ over any finite extension field of $\FF_{p^n}$
actually has coefficients in $\FF_{p^n}$; thus we obtain no new
isomorphisms by making such extensions. Let $\SS_n$ be the group of
automorphisms of this $\Gamma$ over $\FF_{p^n}$; this is the
classical Morava stabilizer group. If we let $\GG_n$ be the group of 
automorphisms of $(\Gamma,\FF_{p^n})$ in $\FGL_n$ (the big Morava stabilizer
group of $\Gamma$), then one easily sees that
$$
\GG_n \cong \SS_n \rtimes \Gal(\FF_{p^n}/\FF_{p}).
$$
Of course, $\GG_n$ acts on $E(\Gamma,\FF_{p^n})$. Also, we note
that the Honda formal group law is defined over $\FF_p$, although
it will not  get its full group of automorphisms until changing
base to $\FF_{p^n}$.

Next we put in the gradings. This requires a paragraph of introduction.
For any commutative ring $R$, the morphism $R[[x]] \to R$ of rings
sending $x$ to $0$ makes $R$ into an $R[[x]]$-module. Let $\Der_R(R[[x]],R)$
denote the $R$-module of continuous $R$-derivations; that is, continuous
$R$-module homomorphisms
$$
\partial : R[[x]] \longrightarrow R
$$
so that
$$
\partial(f(x)g(x)) = \partial(f(x))g(0) + f(0)\partial(g(x)).
$$
If $\partial$ is any derivation, write $\partial(x) = u$; then,
if $f(x) = \sum a_ix^i$,
$$\partial(f(x)) = a_1\partial (x) = a_1u.$$
Thus $\partial$ is determined by $u$, and we write $\partial = \partial_u$.
We then have that\break $\Der_R(R[[x]],R)$ is a free $R$-module of
rank one, generated by any derivation $\partial_u$ so that 
$u$ is a unit in $R$. In the language of schemes, $\partial_u$ is a generator
for the tangent space at $0$ of the formal scheme $\AA^1_R$ over
$\Spec(R)$.

Now consider pairs $(F,u)$ where $F$ is a formal group law over
$R$ and $u$ is a unit in $R$. Thus $F$ defines a smooth one
dimensional commutative formal group scheme over $\Spec(R)$ and
$\partial_u$ is a chosen generator for the tangent space at $0$. A morphism
of pairs
$$
f:(F,u) \longrightarrow (G,v)
$$
is an isomorphism of formal group laws $f:F \to G$ so that
$$
u = f'(0)v.
$$
Note that if $f(x) \in R[[x]]$ is a homomorphism of formal
group laws from $F$ to $G$, and $\partial$ is a derivation at $0$,
then $(f^\ast\partial)(x) = f'(0)\partial(x)$. 
In the context of deformations, we may require that $f$ be a
$\star$-isomorphism.

This suggests the following definition: let $\Gamma$ be a formal group
law of height $n$ over a perfect field $k$ of characteristic $p$,
and let $A$ be a complete local ring. Define $\defo_\Gamma(A)_\ast$
to be equivalence classes of pairs $(G,u)$ where $G$ is a deformation
of $\Gamma$ to $A$ and $u$ is a unit in $A$. 
The equivalence relation is given by 
$\star$-isomorphisms transforming the unit as
in the last paragraph. We now have that there is a natural
isomorphism
$$
\Hom_c(E(\Gamma,k)[u^{\pm 1}],A) \cong \defo_\Gamma(A)_\ast.
$$

We impose a grading by giving an action of the multiplicative
group scheme $\GG_m$ on the scheme $\defo_\Gamma(\cdot)_\ast$ (on the right) 
and thus on $E(\Gamma,k)[u^{\pm 1}]$ (on the left): 
if $v \in A^\times$ is a unit and $(G,u)$ represents an
equivalence class in $\defo_\Gamma(A)_\ast$ define an new
element in $\defo_\Gamma(A)_\ast$ by $(G,v^{-1}u)$. In
the induced grading on $E(\Gamma,k)[u^{\pm 1}]$, one
has $E(\Gamma,k)$ in degree $0$ and $u$ in degree $-2$.

This grading is essentially forced by topological considerations.
See the remarks before Theorem 20 of \cite{Strick} for an explanation.
In particular, it is explained there why $u$ is in degree $-2$ rather
than $2$.

The rest of the section will be devoted to what we need about the
Morava stabilizer group. 
The group $\SS_n$ is the group of units in the endomorphism 
ring $\cO_n$ of the Honda formal group law of height $n$. The ring $\cO_n$ 
can be described as follows (See \cite{Henn} or \cite{Rav}).
One adjoins a noncommuting element $S$ to the
Witt vectors $\WW = W(\FF_{p^n})$ subject to the conditions that
$$
Sa = \phi(a)S\qquad\hbox{and}\qquad S^n = p
$$
where $a \in \WW$ and $\phi:\WW \to \WW$ is the Frobenius. (In terms
of power series, $S$ corresponds to the endomorphism of the formal
group law given by $f(x) = x^p$.) This algebra $\cO_n$ is a free $\WW$-module 
of rank $n$ with generators $1,S,\dots  S^{n-1}$ and is equipped 
with a valuation
$\nu$ extending the standard valuation of $\WW$; since we assume that
$\nu(p) = 1$, we have $\nu(S) = 1/n$. Define a filtration on $\SS_n$
by
$$
F_k\SS_n = \{x \in \SS_n\ |\ \nu(x-1) \geq k\}.
$$
Note that $k$ is a fraction of the form $a/n$ with $a=0,1,2,\dots\ $. We
have
\begin{eqnarray*}
F_{0}\SS_n/F_{1/n}\SS_n &\cong& \FF_{p^n}^\times\ ,
\\
F_{a/n}\SS_n/F_{(a+1)/n}\SS_n &\cong& \FF_{p^n},\qquad a \geq 1
\end{eqnarray*}
and
$$
\SS_n \cong \lim_a \SS_n/F_{a/n}\SS_n.
$$
If we define $S_n = F_{1/n}\SS_n$, then $S_n$ is the $p$-Sylow subgroup
of the profinite group $\SS_n$. Note that the Teichm\"uller elements
$\FF_{p^n}^\times \subseteq \WW^\times \subseteq \cO_n^\times$
define a splitting of the projection $\SS_n \to \FF_{p^n}^\times$ and, hence,
$\SS_n$ is the semi-direct product of $\FF_{p^n}^\times$ and the $p$-Sylow
subgroup.

The action of the Galois group $\Gal(\FF_{p^n}/\FF_p)$ on $\cO_n$ 
is the obvious one: the Galois group is generated by the Frobenius $\phi$ and
$$
\phi(a_0 + a_1S + \dots  + a_{n-1}S^{n-1}) =
\phi(a_0) + \phi(a_1)S + \dots  + \phi(a_{n-1})S^{n-1}.
$$

We are, in this paper, concerned mostly with the case $n=2$ and $p=3$. In
this case, every element of $\SS_2$ can be written as a sum
$$
a + bS, \qquad a,b \in W(\FF_9) = \WW
$$
with $a \not\equiv 0$ mod $3$. The elements of $S_2$ are of the form
$a + bS$ with $a \equiv 1$\break mod $3$.

The following subgroups of $\SS_2$ will be of particular interest
to us. The first two are choices of maximal finite
subgroups.\footnote{The
  first author would like to thank Haynes Miller for
several lengthy and informative discussions about finite subgroups
of the Morava stabilizer group.}
The last one (see \ref{rem0.3})
is a closed subgroup which is, in some sense, complementary to the center.

\begin{sct}\label{rem0.1}
Choose a primitive eighth root of unity $\omega \in \FF_9$. We will write
$\omega$ for the corresponding element in $\WW$ and $\SS_2$. The element
$$
s = -\frac{1}{2}(1+\omega S)
$$
is of order $3${\rm ;} furthermore,
$$
\omega^2s\omega^6 = s^2.
$$
Hence the elements $s$ and $\omega^2$ generate a subgroup of
order $12$ in $\SS_2$ which we label $G_{12}$. As a group, it is abstractly
isomorphic to the unique nontrivial semi-direct product of cyclic groups 
$$
C_3\rtimes C_4.
$$
Any other subgroup of order $12$ in $\SS_2$ is conjugate to $G_{12}$.
In the sequel, when discussing various representations, we will write the
element $\omega^2 \in G_{12}$ as $t$.

We note that the subgroup $G_{12} \subseteq \SS_2$ is a normal
subgroup of a subgroup $G_{24}$ of the larger group $\GG_2$. Indeed,
there is a diagram of short exact sequences of groups
$$
\xymatrix{
{1} \rto & G_{12} \rto \dto_{\subseteq} & G_{24} \rto \dto_{\subseteq}
&\Gal(\FF_9/\FF_3) \rto \dto^{=} &{1}\\
{1} \rto & \SS_2 \rto & \GG_2 \rto & \Gal(\FF_9/\FF_3) \rto &{1}.\\
}
$$
Since the action of the Galois group on $\SS_2$ does not preserve
any choice of $G_{12}$, this is not transparent. In fact, while the lower 
sequence is split the upper sequence is not. More concretely 
we let
$$
\psi = \omega\phi \in \SS_2 \rtimes \Gal(\FF_9/\FF_3) = \GG_2
$$
where $\omega$ is our chosen $8$th root of unity and $\phi$ is the
generator of the Galois group. Then if $s$ and $t$ are the elements
of order $3$ and $4$ in $G_{12}$ chosen above, we easily calculate
that $\psi s = s\psi$, $t\psi = \psi t^3$ and $\psi^2 = t^2$. Thus
the subgroup of $\GG_2$ generated by $G_{12}$ and $\psi$ has order $24$, as
required. Note that the $2$-Sylow subgroup of $G_{24}$ is the
quaternion group $Q_8$ of order $8$ generated by $t$ and $\psi$ and that
indeed 
$$
\xymatrix{
{1} \rto & G_{12} \rto & G_{24} \rto &\Gal(\FF_9/\FF_3) 
\rto& {1}\\
}
$$
is not split.
\end{sct}

\begin{sct}\label{rem0.2}
The second subgroup is the subgroup $SD_{16}$ 
generated by $\omega$ and $\phi$. 
This is the semidirect product 
$$
\FF_9^{\times} \rtimes \ZZ/2\ , 
$$
and it is also known as the semidihedral group of order $16$. 
\end{sct}

\begin{sct}\label{rem0.3}
For the third  subgroup, note that the evident right action of $\SS_n$ on
$\cO_n$ defines a group homomorphism $\SS_n \to \hbox{GL}_n(\WW)$. The
determinant homomorphism $\SS_n \to \WW^\times$ extends to a homomorphism
$$
\GG_n \to \WW^\times\rtimes\Gal(\FF_{p^n}/\FF_p)\ . \\
$$
For example, if $n=2$, this map sends $(a + bS,\phi^e)$, $e\in \{0,1\}$, to
$$
(a\phi(a) - pb\phi(b),\phi^e)
$$
where $\phi$ is the Frobenius. It is simple to check (for all $n$) that the
image of this homomorphism lands in
$$
\ZZ_p^\times\times\Gal(\FF_{p^n}/\FF_p) 
\subseteq \WW^\times \rtimes \Gal(\FF_{p^n}/\FF_p) \ .
$$
If we identify the quotient of 
$\ZZ_p^\times$ by its maximal finite subgroup
with $\ZZ_p$,  
we get a ``reduced determinant'' homomorphism
$$
\GG_n \to \ZZ_p\ .
$$
Let $\GG^1_n$ be the kernel of this map and $\SS_n^1$ resp. 
$S_n^1$ be the kernel of its restriction to 
$\SS_n$ resp. $S_n$. 
In particular, any finite subgroup of $\GG_n$ is a subgroup of $\GG_n^1$. 
One also easily checks that the center of $\GG_n$ is $\ZZ_p^\times
\subseteq \WW^\times \subseteq \SS_n$ and that the composite
$$
\ZZ_p^\times  \to \GG_n \to \ZZ_p^\times
$$
sends $a$ to $a^n$. Thus, if $p$ does not divide $n$, we have
$$
\GG_n \cong \ZZ_p \times \GG_n^1.
$$ 
\end{sct}

\section{The $K(n)$-local category and the Lubin-Tate theories $E_n$}

The purpose of this section is to collect together the information
we need about the $K(n)$-local category and the role of the functor
$(E_n)_\ast(\cdot)$ in governing this category. {\it But
attention\/}! --- $(E_n)_\ast X$ is {\it not} the homology
of $X$ defined by the spectrum $E_n$, but a 
completion thereof; see Definition \ref{defn1.1} below.

Most of the information in this section is collected from
\cite{Dev}, \cite{DH}, and \cite{HS}.

Fix a prime $p$ and let $K(n)$, $1 \leq n < \infty$, denote the
$n$-th Morava $K$-theory spectrum. Then $K(n)_\ast \cong \FF_p[v_n^{\pm 1}]$
where the degree of $v_n$ is $2(p^n-1)$. This is a complex
oriented theory and the formal group law over $K(n)_\ast$ is
of height $n$. As is customary, we specify that the formal group law over 
$K(n)_\ast$ is the graded variant of the Honda formal group law;
thus, the $p$-series is
$$
[p](x) = v_n x^{p^n}.
$$

Following Hovey and Strickland, we will write $\cK_n$ for
the category of\break $K(n)$-local spectra. We will write $L_{K(n)}$
for the localization functor from spectra to $\cK_n$. 

Next let $K_n$ be the extension of $K(n)$ with $(K_n)_\ast \cong
\FF_{p^n}[u^{\pm 1}]$ with the degree of $u = -2$. The
inclusion $K(n)_\ast \subseteq (K_n)_\ast$ sends $v_n$ to
$u^{-(p^n-1)}$. There is a natural isomorphism of homology theories
$$
(K_n)_\ast \otimes_{K(n)_\ast} K(n)_\ast X \mathop{\longrightarrow}^{\cong}
(K_n)_\ast X
$$
and $K(n)_\ast \to (K_n)_\ast$ is a faithfully flat extension; thus
the two theories have the same local categories and weakly
equivalent localization functors.

If we write $F$ for the graded formal group law over $K(n)_\ast$
we can extend $F$ to a formal group law over $(K_n)_\ast$ and define
a formal group law $\Gamma$ over $\FF_{p^n} = (K_n)_0$ by
$$
x+_{\Gamma}y = \Gamma(x,y) = u^{-1}F(ux,uy) = u^{-1}(ux +_F uy).
$$
Then $F$ is chosen so that $\Gamma$ is the Honda formal group law.
\smallbreak 

We note that -- as in \cite{DH} -- there is a choice of the 
universal deformation $G$ of $\Gamma$ such that the $p$-series 
of the associated graded formal group law $G_0$ over 
$E(\Gamma,\FF_{p^n})[u^{\pm 1}]$ satisfies 
$$
[p](x) = v_0x +_{G_0} v_1x^p +_{G_0} v_2x^{p^2} +_{G_0} \dots 
$$
with $v_0=p$ and
$$
v_k = \left\{ \begin{array}{ll}
u^{1-p^k}u_k & 0 < k < n;\\
u^{1-p^n} & k = n;\\
0 & k > n.\\
\end{array}\right.
$$

This shows that the functor 
$X\mapsto (E_n)_\ast\otimes_{BP_*}BP_*X$ (where $(E_n)_{\ast}$ is considered 
a $BP_{\ast}$-module via the evident ring homomorphism) 
is a homology theory which is represented by a spectrum $E_n$ with 
coefficients 
$$
\pi_*(E_n)\cong E(\Gamma,\FF_{p^n})[u^{\pm 1}]\cong
\WW[[u_1,\dots ,u_{n-1}]][u^{\pm 1}]\ .
$$
The inclusion of the subring $E(n)_*=\ZZ_{(p)}[v_1,\dots ,v_{n-1},v_n^{\pm 1}]$ 
into $(E_n)_\ast$ is again faithfully
flat; thus, these two theories have the same local categories. We
write $\cL_n$ for the category of $E(n)$-local spectra and $L_n$
for the localization functor from spectra to $\cL_n$.

The reader will have noticed that we have avoided using the
expression $(E_n)_\ast X$; we now explain what we mean by this.
The $K(n)$-local category $\cK_n$ has internal smash products and  
(arbitrary) wedges given by
$$
X \wedge_{\cK_n} Y = L_{K(n)}(X \wedge Y)
$$
and 
$$
\mathop{\bigvee}_{\cK_n} X_\alpha = L_{K(n)}(\mathop{\bigvee} X_\alpha)\ .
$$

In making such definitions, we assume we are working in some suitable
model category of spectra, and that we are taking the smash
product between cofibrant spectra; that is, we are working with the
derived smash product. The issues here are troublesome, but
well understood, and we will not dwell on these points.
See \cite{EKMM} or \cite{HSS}. 
If we work in our suitable categories of spectra 
the functor $Y \mapsto X \wedge_{\cK_n} Y$ has
a right adjoint $Z \mapsto F(X,Z)$.

We define a version of $(E_n)_\ast(\cdot)$ intrinsic to $\cK_n$
as follows.

\begin{defn}\label{defn1.1} Let $X$ be a spectrum. Then we define
$(E_n)_\ast X$ by the equation
$$
(E_n)_\ast X = \pi_\ast L_{K(n)}(E_n \wedge X).
$$
\end{defn}

We remark immediately that $(E_n)_\ast (\cdot)$ is not a homology
theory in the usual sense; 
for example, it will not send arbitrary wedges
to sums of abelian groups. However, it is tractable, as we now
explain. First note that $E_n$ itself is $K(n)$-local; indeed,
Lemma 5.2 of \cite{HS} demonstrates that $E_n$ is a finite wedge
of spectra of the form $L_{K(n)}E(n)$. Therefore if $X$ is a finite
$CW$ spectrum, then $E_n \wedge X$ is already in $\cK_n$, so
\begin{equation}\label{eq1.1}
(E_n)_\ast X = \pi_\ast (E_n \wedge X).
\end{equation}

Let $I = (i_0,\dots ,i_{n-1})$ be a sequence of positive integers
and let
$$
\mm^I = (p^{i_0},u_1^{i_1},\dots ,u_{n-1}^{i_{n-1}}) \subseteq \mm \subseteq
(E_n)_\ast
$$
where $\mm=(p,u_1,\dots ,u_{n-1})$ is the maximal ideal in $E_\ast$. 
These form a system of ideals in $(E_n)_\ast$ and produce a filtered
diagram of rings $\{(E_n)_\ast/\mm^I\}$; furthermore
$$
(E_n)_\ast = \lim_I (E_n)_\ast/\mm^I.
$$
There is a cofinal diagram $\{(E_n)_\ast/\mm^J\}$ which can be realized
as a diagram of spectra in the following sense: using nilpotence
technology, one can produce a diagram
of finite spectra $\{M_J\}$ and an isomorphism
$$
\{(E_n)_\ast M_J\} \cong \{(E_n)_\ast/\mm^J\}
$$
as diagrams. See \S 4 of \cite{HS}. Here $(E_n)_\ast M_J =
\pi_\ast E_n \wedge M_J = \pi_\ast L_K(n)(E_n \wedge M_J)$.
The importance of this diagram is that (see \cite[Prop.~7.10]{HS}) 
for each spectrum X 
\begin{equation}\label{eq1.2} 
L_{K(n)}X \simeq \holim_J M_J \wedge L_nX.
\end{equation}
This has the following consequence, immediate from Definition \ref{defn1.1}:
there is a short exact sequence
$$
0 \to \lim{}^1 (E_n)_{k+1}(X \wedge M_J) \to (E_n)_k X \to
\lim (E_n)_{k}(X \wedge M_J) \to 0.
$$
This suggests $(E_n)_\ast X$ is closely related to some
completion of $\pi_\ast (E_n \wedge X)$ and this is nearly the
case. The details are spelled out in Section  8 of \cite{HS}, but we
will not  need the full generality there. In fact, all of the
spectra we consider here will satisfy the hypotheses of Proposition
\ref{prop1.2} below.

If $M$ is an $(E_n)_\ast$-module, let $M_\mm^{\wedge}$ denote the completion
of $M$ with respect to the maximal ideal of $(E_n)_\ast$. 
A module of the form 
$$
(\bigoplus_\alpha \Sigma^{k_{\alpha}}(E_n)_\ast)^{\wedge}_\mm
$$
will be called {\it pro-free}.

\begin{prop}\label{prop1.2} If $X$ is a spectrum so that $K(n)_\ast X$
is concentrated in even degrees{\rm ,} then
$$
(E_n)_\ast X \cong \pi_\ast (E_n \wedge X)^{\wedge}_\mm
$$
and $(E_n)_\ast X$ is pro-free as an $(E_n)_\ast$\/{\rm -}\/module.
\end{prop}

See Proposition 8.4 of \cite{HS}.
 \smallbreak 

As with anything like a flat homology theory, the object 
$(E_n)_\ast X$
is a comodule over some sort of Hopf algebroid of co-operations;
it is our next project to describe this structure. In particular,
this brings us to the role of the Morava stabilizer group. We begin
by identifying $(E_n)_\ast E_n$.

Let $\GG_n$ be the (big) Morava
stabilizer group of $\Gamma$, the Honda formal group law of height
$n$ over $\FF_{p^n}$. For the purposes of this paper, a {\it Morava module}
is a complete $(E_n)_\ast$-module $M$ equipped with a continuous 
$\GG_n$-action subject to the following compatibility condition: if 
$g \in \GG_n$, $a \in (E_n)_\ast$ and $x \in M$,
then
\begin{equation}\label{twistedaction}
g(ax) = g(a)g(x)\ .
\end{equation}
For example, if $X$ is any spectrum with $K(n)_\ast X$ concentrated in
even degrees, then $(E_n)_\ast X$ is a complete $(E_n)_\ast$-module
(by Proposition \ref{prop1.2}) and the action of  $\GG_n$ on $E_n$ defines a
continuous action of $\GG_n$ on $(E_n)_\ast X$. This is a prototypical
Morava module.

Now let $M$ be a Morava module and let
$$
\Hom^c(\GG_n,M)
$$
be the abelian group of continuous maps from $\GG_n$ to $M$ 
where the topology on $M$ is defined via the ideal $\mm$. Then
\begin{equation}\label{eq1.4}
\Hom^c(\GG_n,M) \cong \lim{}_i\colim_{k}\map(\GG_n/U_k,M/\mm^iM)
\end{equation}
where $U_k$ runs over any system of open subgroups of $\GG_n$ with
$\bigcap_k U_k = \{e\}$.
To give $\Hom^c(\GG_n,M)$ a structure of an $(E_n)_\ast$-module
let $\phi:\GG_n \to M$ be continuous and $a \in (E_n)_\ast$.
The we define $a\phi$ by the formula
\begin{equation}\label{module}
(a\phi)(x) = a\phi(x)\ .
\end{equation}
There also is a continuous action of
$\GG_n$ on $\Hom^c(\GG_n,M)$: if $g \in \GG_n$ and $\phi:\GG_n \to M$
is continuous, then one defines $g\phi:\GG_n \to M$ by the formula
\begin{equation}\label{action}
(g\phi)(x) = g\phi(g^{-1}x)\ .
\end{equation}
With this action, and the action of $(E_n)_\ast$ defined in 
(\ref{module}),  the formula of (\ref{twistedaction}) holds.
Because $M$ is complete (\ref{eq1.4}) shows that $\Hom^c(\GG_n,M)$
is complete.

\begin{rem}\label{right-adjoint} With the Morava module structure
defined by equations \ref{module} and \ref{action}, the functor
$M \to \Hom^c(\GG_n,M)$ has the following universal property.
If $N$ and $M$ are Morava modules and $f:N \to M$ is a morphism
of continuous $(E_n)_\ast$ modules, then there is an induced
morphism
\begin{align*}
N &\longrightarrow \Hom^c(\GG_n,M)\\
\alpha &\mapsto \phi_\alpha 
\end{align*}
with $\phi_\alpha(x) = xf(x^{-1}\alpha)$. This yields a natural isomorphism
$$
\Hom_{(E_n)_\ast}(N,M) = \Hom_{\rm Morava}(N,\Hom^c(\GG_n,M))
$$
from continuous $(E_n)_\ast$ module homomorphisms to morphisms
of Morava modules.

There is a {\it different}, but isomorphic natural Morava module structure
on $\Hom^c(\GG_n,-)$ so that this functor becomes a true right adjoint
of the forget functor from Morava modules to continuous
$(E_n)_\ast$-modules. However, we will not need this module structure at 
any point and we supress it to avoid confusion.
\end{rem}

For example, if $X$ is a spectrum such that $(E_n)_\ast X$ is 
$(E_n)_\ast$-complete, the $\GG_n$-action on $(E_n)_\ast X$ is encoded
by the map
$$
(E_n)_\ast X \to \Hom^c(\GG_n,(E_n)_\ast X)
$$
adjoint (in the sense of the previous remark) to the identity. 

The next result says that this is essentially all the stucture that
$(E_n)_\ast X$ supports. For any spectrum $X$, $\GG_n$ acts on 
$$
(E_n)_\ast(E_n \wedge X) = \pi_\ast L_{K(n)}(E_n\wedge E_n \wedge X)
$$
by operating in the left factor of $E_n$. 
The multiplication $E_n \wedge E_n \to E_n$ defines
a morphism of $(E_n)_\ast$-modules
$$
(E_n)_\ast (E_n \wedge X) \to (E_n)_\ast X
$$ 
and by composing we obtain a map 
$$
\phi:(E_n)_\ast (E_n \wedge X) \to \Hom^c(\GG_n,(E_n)_\ast (E_n\wedge X)) 
\to \Hom^c(\GG_n,(E_n)_\ast X)\ .
$$
If $(E_n)_\ast X$ is complete, this is a morphism of Morava modules. 

We now record:

\begin{prop}\label{prop1.3} For any 
cellular spectrum $X$ with $(K_n)_\ast X$
concentrated in even degrees the morphism
$$
\phi:(E_n)_\ast (E_n \wedge X) \to \Hom^c(\GG_n,(E_n)_\ast X)
$$
is an isomorphism of Morava modules.
\end{prop}

\Proof  See \cite{DH1} and \cite{Strick} for the case $X = S^0$.
The general case
follows in the usual manner. First, it's true for finite spectra by
a five lemma argument. For this one needs to know that the functor
$$
M \mapsto \Hom^c(\GG_n,M)
$$
is exact on finitely generated $(E_n)_*$-complete modules. This follows
from (\ref{eq1.4}).  
Then one argues the general case, by noting
first that by taking colimits 
over finite cellular subspectra
$$
\phi:(E_n)_\ast (E_n \wedge M_J \wedge X) \to
\Hom^c(\GG_n,(E_n)_\ast (M_J \wedge X))
$$
is an isomorphism for any $J$ and any $X$. Note that
$E_n\wedge M_J\wedge X$ is $K(n)$-local for any $X$; therefore,
$L_{K(n)}$ commutes with the homotopy colimits in question. 
Finally the hypothesis on $X$ implies 
$$
(E_n)_\ast(E_n \wedge X) \cong \lim (E_n)_\ast (E_n \wedge M_J \wedge X). 
$$
and thus we can conclude the result
by taking limits with respect to $J$. 
\Endproof\vskip4pt

We next turn to the results of Devinatz and Hopkins (\cite{DH1}) on homotopy
fixed point spectra. 
Let $\cO_{\GG_n}$ be
the orbit category of $\GG_n$. Thus an object in $\cO_{\GG_n}$ is
an orbit $\GG_n/H$ where $H$ is a closed subgroup and the morphisms
are continuous $\GG_n$-maps. Then Devinatz and Hopkins have defined
a functor
$$
\cO_{\GG_n}^{op} \to \cK
$$
sending $\GG_n/H$ to a $K(n)$-local spectrum $E_n^{hH}$. If $H$ is 
finite, then $E_n^{hH}$ is the usual homotopy fixed point spectrum
defined by the action of $H \subseteq \GG_n$. By the results
of \cite{DH1}, the morphism $\phi$ of Proposition \ref{prop1.3}
restricts to an isomorphism (for any closed $H$) 
\begin{equation}\label{fixediso}
(E_n)_\ast E_n^{hH} \mathop{\longrightarrow}^{\cong}
\Hom^c(\GG_n/H,(E_n)_\ast).
\end{equation}

We would now like to write down a result about the function
spectra $F((E_n)^{hH},E_n)$. First, some notation. If
$E$ is a spectrum and $X=\lim_{i}X_{i}$
is an inverse limit of a sequence of finite sets $X_i$ then define
$$
E[[X]] = \holim_{i}E\wedge (X_{i})_{+}.
$$

\begin{prop}\label{functioniso-1} Let $H$ be a closed subgroup of
$\GG_n$. Then there is a natural weak equivalence
$$
\xymatrix{
E_n[[\GG_n/H]] \rto^\simeq & F((E_n)^{hH},E_n).\\
}
$$
\end{prop}

\Proof  First let $U$ be an open subgroup of $\GG_n$. 
Functoriality of the homotopy fixed point spectra
construction of \cite{DH1}
gives us a map $E_n^{hU}\wedge \GG_n/U_{+}\to E_n$
where as usual $\GG_n/U_+$ denotes $\GG_n/U$ with a disjoint base point
added. Together with the product on $E_n$ we obtain a map
\begin{equation}\label{multiplication}
E_n\wedge E_n^{hU}\wedge \GG_n/U_{+}\to E_n\wedge E_n\to E_n
\end{equation}
whose adjoint induces an equivalence
\begin{equation}\label{point-set-fixediso}
L_{K(n)}(E_n\wedge E_n^{hU})\to \prod_{\GG_n/U}E_n
\end{equation}
realizing the isomorphism of (\ref{fixediso}) above. Note that this
is a map of $E_n$-module spectra. Let $F_{E_n}(-,E_n)$
be the function spectra in the category of left $E_n$-module
spectra. (See \cite{EKMM} for details.) If we apply
$F_{E_n}(-,E_n)$ to the equivalence of (\ref{point-set-fixediso}) we
obtain an equivalence of $E_n$-module spectra
$$
F_{E_n}(\prod_{\GG_n/U} E_n, E_n) \to F_{E_n}(E_n \wedge E_n^{hU}, E_n).
$$
This equivalence can then be written as
\begin{equation}\label{open-to-the-rescue}
E_n \wedge (\GG_n/U)_+ \to F(E_n^{hU}, E_n);
\end{equation}
furthermore, an easy calculation shows that this map is adjoint to
the map of (\ref{multiplication}). 

More generally, let $H$ be any closed subgroup of $\GG_n$. Then there exists
a decreasing sequence $U_i$ of open subgroups $U_i$ with $H=\bigcap_iU_i$
and by \cite{DH1} we have
$$
E_n^{hH}\simeq L_{K(n)}\hocolim_iE_n^{hU_i} \ .
$$
Thus, the equivalence of (\ref{open-to-the-rescue}) and by passing to
the limit we obtain the desired equivalence.
\Endproof\vskip4pt 

Now note that if $X$ is a profinite set with continuous
$H$-action and if $E$ is an $H$-spectrum then $E[[X]]$ is an
$H$-spectrum via the diagonal action. It is this action which
is used in the following result. 

\begin{prop}\label{prop1.6} {\rm 1)} Let $H_1$ be a closed subgroup and $H_2$ a
finite subgroup of $\GG_n$. Then there is a natural equivalence
$$
\xymatrix{
E_n[[\GG_n/H_1]]^{hH_2}\rto^\simeq& F(E_n^{hH_1},E_n^{hH_2})\ .\\
}
$$

{\rm 2)} If $H_1$ is also an open subgroup then there is a natural decomposition
$$
E_n[[\GG_n/H_1]]^{hH_2} \simeq \prod_{H_2\backslash \GG_n/H_1} E_n^{hH_x}
$$
where $H_x = H_2 \cap xH_1 x^{-1}$ is the isotropy subgroup of the coset 
$xH_1$ and $H_2\backslash \GG_n/H_1$ is the \/{\rm (}\/finite\/{\rm )}\/ set of double cosets.
\smallbreak
{\rm 3)} If $H_1$ is a closed subgroup and $H_1=\bigcap_iU_i$ for
a decreasing sequence of open subgroups $U_i$
then
$$
F(E_n^{hH_1},E_n^{hH_2}) \simeq \holim_iE_n[[\GG_n/U_i]]^{hH_2}\simeq
\holim_i\prod_{H_2\backslash \GG_n/U_i} E_n^{hH_{x,i}} 
$$
where $H_{x,i} = H_2 \cap xU_ix^{-1}$ is{\rm ,} as before{\rm ,}
the isotropy subgroup of the coset $xU_i$.\footnote{We are grateful to P. Symonds for pointing out  
that the naive generalization of the second 
statement does not hold for a general closed subgroup.}
\end{prop}

\Proof  The first statement follows from Proposition
\ref{functioniso-1} by passing to homotopy fixed point spectra 
with respect to $H_2$ and the second statement is then an 
immediate consequence of the first.  For the third statement 
we write $\GG_n/H_1 = \lim_i \GG_n/U_i$ 
and pass to the homotopy inverse limit. 
\Endproof\vskip4pt 

We will be interested in the $E_n$-Hurewicz homomorphism
$$
\pi_0 F(E_n^{hH_1},E_n^{hH_2}) \to
\Hom_{(E_n)_\ast E_n}((E_n)_\ast E_n^{hH_1},(E_n)_\ast E_n^{hH_2})
$$
where $\Hom_{(E_n)_\ast E_n}$ denotes morphisms in the category
of Morava modules. Let
$$
(E_n)_\ast[[\GG_n]] = \lim_i(E_n)_\ast[\GG_n/U_i]
$$
denote the completed group ring and give this the structure
of a Morava module by letting $\GG_n$ act diagonally. 

\begin{prop}\label{prop1.7} Let $H_1$ and $H_2$ be closed subgroups
of $\GG_n$ and suppose that $H_2$ is finite. Then there is an isomorphism 
$$
\big((E_n)_\ast[[\GG_n/H_1]]\big)^{H_2} \mathop{\longrightarrow}^{\cong} 
\Hom_{(E_n)_\ast E_n}((E_n)_\ast E_n^{hH_1},(E_n)_\ast E_n^{hH_2})
$$
such that the following diagram commutes
$$
\xymatrix{
\pi_*E_n[[\GG_n/H_1]]^{hH_2} \rto \dto_{\cong} & 
\big((E_n)_*[[\GG_n/H_1]]\big)^{H_2} \dto^\cong\\
\pi_*F(E_n^{hH_1},E_n^{hH_2}) \rto & 
\Hom_{(E_n)_\ast E_n}((E_n)_\ast 
E_n^{hH_1},(E_n)_\ast E_n^{hH_2}) 
}
$$
where the top horizontal 
map is the edge homomorphism in the homotopy fixed point
spectral sequence{\rm ,} the left\/{\rm -}\/hand vertical map is induced by the 
equivalence of Proposition {\rm \ref{prop1.6}} and the bottom horizontal map 
is the $E_n$\/{\rm -}\/Hurewicz homomorphism.
\end{prop}

\Proof  First we assume that $H_2$ is the trivial subgroup
and $H_1$ is open, so that $\GG_n/H_1$ is finite. Then there is
an isomorphism
$$
(E_n)_\ast[[\GG_n/H_1]] \to \Hom_{(E_n)_\ast} (\Hom^c (\GG_n/H_1, (E_n)_\ast),
(E_n)_\ast)
$$
which is the unique linear map which sends a coset to evaluation
at that coset. Applying Remark \ref{right-adjoint} we obtain an isomorphism
of Morava modules
$$
(E_n)_\ast[[\GG_n/H_1]] \to \Hom_{(E_n)_\ast} (\Hom^c (\GG_n/H_1, (E_n)_\ast),
\Hom^c (\GG_n,(E_n)_\ast)).
$$
This isomorphism can be extended to a general closed subgroup $H_1$ by
writing $H_1$ as the intersection of a nested sequence of open
subgroups (as in the proof of Proposition \ref{functioniso-1}) and
taking limits.
Then we use the isomorphisms of (\ref{fixediso}) to identify 
$(E_n)_*E_n^{hH_i}$ with $\Hom^c(\GG_n/H_i,(E_n)_*)$. This defines
the isomorphism we need, and it is straightforward to see that the
diagram commutes. To end the proof, note that the case of a general
finite subgroup $H_2$ follows by passing to $H_2$-invariants. 
\hfill\qed

\section{The homotopy groups of $E_2^{hF}$ at $p=3$}

To construct our tower we are going to need some information about
$\pi_{\ast}E^{hF}_2$ for various finite subgroups of the stabilizer
group $\GG_2$. Much of what we say here can be recovered from various places
in the literature (for example, \cite{GS}, \cite{Nave}, or \cite{GHM}) and
the point of view and proofs expressed are certainly those of Mike Hopkins. 
What we add here to the discussion in \cite{GHM} is 
that we pay careful attention to the Galois group. In particular we treat 
the case of the finite group $G_{24}$.

Recall that we are working at the prime 3. We will write $E$ for
$E_2$, so that we may write $E_\ast$ for $(E_2)_\ast$.

In Remark \ref{rem0.1} we defined a subgroup
$$
G_{24} \subseteq \GG_2 = \SS_2 \rtimes \Gal(\FF_9/\FF_3)
$$
generated by elements $s$, $t$ and $\psi$ of orders $3$, $4$ and $4$
respectively. The cyclic subgroup $C_3$ generated by $s$ is normal,
and the subgroup $Q_8$ generated by $t$ and $\psi$ is the quaternion
group of order $8$.

The first results are algebraic in nature; they give a nice
presentation of $E_\ast$ as a $G_{24}$-algebra. First we define an
action of $G_{24}$ on $\WW = W(\FF_9)$ by the formulas:
\begin{equation}\label{main-representation}
s(a) = a\qquad t(a) = \omega^2a\qquad \psi(a) = \omega\phi(a)
\end{equation}
where $\phi$ is the Frobenius. Note the action factors through
$G_{24}/C_3 \cong Q_8$. Restricted to the subgroup $G_{12} = \SS_2 \cap
G_{24}$ this action is $\WW$-linear, but over $G_{24}$ it is
simply linear over $\ZZ_3$. Let $\Chi$ denote the resulting
$G_{24}$-representation and $\Chi'$ its restriction to $Q_8$. 

This representation is a module over a twisted version of the
group ring $\WW[G_{24}]$. The projection
$$
G_{24} \longrightarrow {\rm Gal}(\FF_9/\FF_3)
$$
defines an action\footnote{This action
is different from that of the representation defined by the formulas of
\ref{main-representation}.} of $G_{24}$ on $\WW$ and we use this
action to twist the multiplication in $\WW[G_{24}]$.
We should really write $\WW_{\phi}[G_{24}]$ for this twisted group
ring, but we forebear, so as to not clutter notation. Note
that $\WW[Q_8]$ has a similar twisting, but $\WW[G_{12}]$ is the
ordinary group ring.

Define a $G_{24}$-module $\rho$ by the short exact sequence
\begin{equation}\label{eq100.1}
0 \to \Chi \to \WW[G_{24}] \otimes_{\WW[Q_8]} \Chi' \to \rho \to 0
\end{equation}
where the first map takes a generator $e$ of $\Chi$ to
$$
(1+s+s^2)e \in \WW[G_{24}] \otimes_{\WW[Q_8]} \Chi'.
$$

\begin{lem}\label{lem100.1} There is a morphism of $G_{24}$-modules
$$
\rho \longrightarrow E_{-2}
$$
so that the induced map
$$
\FF_9 \otimes_{\WW} \rho \to E_0/(3,u_1^2) \otimes_{E_0} E_{-2}
$$
is an isomorphism. Furthermore{\rm ,} this isomorphism sends the
generator $e$ of $\rho$ to an invertible element in $E_\ast$.
\end{lem}

\Proof  We need to know a bit about the action of $\GG_2$ on $E_\ast$. 
The
relevant formulas have been worked out by Devinatz and Hopkins. Let
$\mathfrak{m} \subseteq E_0$ be the maximal ideal and $a + bS \in \SS_2$. Then
Proposition 3.3 and Lemma 4.9 of \cite{DH} together imply that,
modulo $\mathfrak{m}^2E_{-2}$
\begin{align}\label{devhop}
(a+bS)u &\equiv au + \phi(b)uu_1 \\
(a+bS)uu_1 &\equiv 3bu + \phi(a)uu_1 \ .
\end{align}
In some cases we can be more specific. For example, if $\alpha \in
\FF_9^\times \subseteq \mathbb{W}^\times \subseteq
\GG_2$, then the induced map of rings
$$
\alpha_\ast:E_\ast \to E_\ast
$$
is the $\WW$-algebra map defined by the formulas
\begin{equation}\label{action-of-fstar}
\alpha_\ast (u) = \alpha u\qquad\hbox{and}\qquad \alpha_\ast (uu_1) =
\alpha^3uu_1\ .
\end{equation}
Finally, since the Honda formal group is defined over $\FF_3$ the action
of the Frobenius on $E_\ast = \WW [[u_1]][u^{\pm 1}]$ is
simply extended from the action on $\WW$. Thus we have
\begin{equation}\label{devhop1}
\psi(x) = \omega_\ast\phi(x)
\end{equation}
for all $x\in E_2$.

The formulas (\ref{devhop}) up to (\ref{devhop1}) 
imply that $E_0/(3,u_1^2) \otimes_{E_0} E_{-2}$
is isomorphic to $\FF_9 \otimes_{\WW} \rho$ as a $G_{24}$-module and, further,
that we can choose as a generator the residue class of
$u$. In \cite{GHM} (following \cite{Nave}, who learned it from Hopkins) we
found a class $y \in E_{-2}$ so that
\begin{equation}\label{eq100.3}
y \equiv \omega u  \qquad\hbox{mod}\qquad  (3,u_1) 
\end{equation}
and so that
$$
(1 + s + s^2)y = 0 . 
$$
This element might not yet have the correct invariance property with
respect to $\psi$; to correct this, we average and set 
\begin{eqnarray*}
x& =& \frac{1}{8}(y + \omega^{-2}t_\ast(y) + \omega^{-4}(t^2)_\ast(y)
+ \omega^{-6}(t^3)_\ast(y) \\
&& +\ \omega^{-1}\psi_\ast(y) 
 + \omega^{-7}(\psi t)_\ast(y) + \omega^{-5}(\psi t^2)_\ast (y)
+ \omega^{-3}(\psi t^3)_\ast (y)).
\end{eqnarray*}
We can now send the generator of $\rho$ to $x$. Note also that the
formulas (\ref{devhop}) up to  (3.7) 
imply that 
\vglue12pt
\hfill $
\displaystyle{x \equiv \frac{1}{8}(\omega u+\omega^3u)\ \text{modulo} \ (3,u_1^2)\ .}
$ 
\Endproof\vskip12pt

We now make a construction. The morphism of $G_{24}$-modules constructed
in this last lemma defines a morphism of $\WW$-algebras
$$
S(\rho)= S_\WW(\rho) \longrightarrow E_\ast
$$
sending the generator $e$ of $\rho$ to an invertible element in $E_{-2}$.
The symmetric algebra is over $\WW$ and the map is a map of $\WW$-algebras.
The group $G_{24}$ acts through $\ZZ_3$-algebra maps, and the subgroup
$G_{12}$ acts through $\WW$-algebra maps. If $a \in \WW$ is a multiple of
the unit, then $\psi(a) = \phi(a)$.

Let
\begin{equation}\label{eq100.4}
N = \prod_{g\in G_{12}} ge \in S(\rho);
\end{equation}
then $N$ is invariant by $G_{12}$ and $\psi(N)=-N$ so that we get a morphism of graded
$G_{24}$-algebras 
$$
S(\rho)[N^{-1}] \longrightarrow E_\ast
$$
(where the grading on the source is determined by putting $\rho$ 
in degree $-2$). 
Inverting $N$ inverts $e$, but in an invariant manner.  This map is not
yet an isomorphism, but it is an inclusion onto a dense subring.
The following result is elementary (cf. Proposition
2 of \cite{GHM}):

\begin{lem}\label{lem100.2} Let $I = S(\rho)[N^{-1}] \cap \mm$. Then
completion at the ideal $I$ defines an
isomorphism of $G_{24}$-algebras
$$
S(\rho)[N^{-1}]^\wedge_I \cong E_\ast.
$$
\end{lem}

\def\tran{{{\mathrm{tr}}}}

Thus the input for the calculation of the $E_2$-term 
$H^\ast (G_{24},E_\ast)$ of the homotopy fixed point spectral 
sequence associated to $E_2^{hG_{24}}$ will
be discrete. Indeed, let $A = S(\rho)[N^{-1}]$. Then the essential
calculation is that of $H^\ast (G_{24},A)$. For this 
we begin with the following.
For any finite group $G$ and any $G$ module~$M$, let
$$
\tran_G = \tran: M \longrightarrow M^G = H^0(G,M)
$$
be the transfer: $\tran(x) = \sum_{g\in G} gx$. In the following result,
an element listed as being in bidegree $(s,t)$ is in $H^s(G,A_t)$.

If $e \in \rho$ is the generator, define $d \in A$ to be the 
multiplicative norm with respect to the cyclic group $C_3$ generated by $s$: 
$d = s^2(e)s(e)e$. By construction $d$ is 
invariant with respect to $C_3$. 

\begin{lem}\label{lem100.4}Let $C_3 \subseteq G_{12}$ be the normal
subgroup of order three. Then there is an exact sequence
$$
A \mathop{\longrightarrow}^{\tran} H^\ast (C_3,A)
\to \FF_9[a,b,d^{\pm 1}]/(a^2) \to 0
$$
where $a$ has bidegree $(1,-2)${\rm ,} $b$ has bidegree $(2,0)$ and $d$ has
bidegree $(0,-6)$. Furthermore the action of $t$ and $\psi$ is described
by the formulas
$$
t(a) = -\omega^2a \qquad t(b) = -b \qquad t(d)
= \omega^6d
$$
and
$$
\psi(a) = \omega a \qquad \psi(b) = b \qquad \psi(d) = \omega^3d \ .
$$
\end{lem}

\Proof  This is the same argument as in Lemma 3 of \cite{GHM},
although here we keep track of the Frobenius.

Let $F$ be the $G_{24}$-module 
$\WW[G_{24}] \otimes_{\WW[Q_8]} \Chi'$; thus equation \ref{eq100.1}
gives a short exact sequence of $G_{24}$-modules
\begin{equation}\label{eq100.5}
0 \to S(F)\otimes \Chi \to S(F) \to S(\rho) \to 0\ .
\end{equation}
In the first term, we set the degree of $\Chi$ to be $-2$ in order
to make this an exact sequence of graded modules. We use the resulting long
exact sequence for computations. We may choose $\WW$-generators of
$F$ labelled $x_1$, $x_2$, and $x_3$ so that if $s$ is the chosen
element of order 3 in $G_{24}$, then $s(x_1) = x_2$ and
$s(x_2) = x_3$. Furthermore, we can choose $x_1$ so that it maps to
the generator $e$ of $\rho$ and is invariant under the action
of the Frobenius. Then we have
$$
S(F) = \WW[x_1,x_2,x_3]
$$
with the $x_i$ in degree $-2$. Under the action of $C_3$ the orbit
of a monomial in $\WW[x_1,x_2,x_3]$ has three elements unless that
monomial is a power of $\sigma_3 = x_1x_2x_3$ -- which, of course,
maps to $d$. Thus, we have a short exact sequence
$$
S(F) \mathop{\longrightarrow}^{\tran} H^\ast (C_3,S(F)) \to
\FF_9[b,d] \to 0
$$
where $b$ has bidegree $(2,0)$ and $d$ has bidegree $(0,-6)$.
Here $b \in H^2(C_3,\ZZ_3) \subseteq H^2(C_3,\WW)$ is a generator and
$\WW \subseteq S(F)$ is the submodule generated by the algebra
unit. Note that the action of $t$ is described by
$$
t(d) = \omega^6 d \qquad\hbox{and}\qquad t(b) = -b\ .
$$
The last is because the element $t$ acts nontrivially on
the subgroup $C_3 \subseteq G_{24}$ and hence on $H^2(C_3,\WW)$.
Similarly, since the action of the Frobenius on $d$ is trivial and
$\psi$ acts trivially on $C_3$, we have
$$
\psi (d) = \omega^3d\qquad\hbox{and}\qquad \psi(b)=b\ .
$$
The short exact sequence (\ref{eq100.5}) 
and the fact that $H^1(C_3,S(F)) = 0$ now imply 
that there is an exact sequence
$$
S(\rho) \mathop{\longrightarrow}^{\tran} H^\ast (C_3,S(\rho)) \to
\FF_9[a,b,d]/(a^2) \to 0\ .
$$
The element $a$ maps to
$$
b \in H^2 (C_3,S_0(F) \otimes \Chi) = H^2 (C_3,\Chi)
$$
under the boundary map (which is an isomorphism)
$$
H^1 (C_3,\rho) = H^1 (C_3,S_1(\rho)) \to H^2 (C_3,\Chi);
$$
thus $a$ has bidegree $(1,-2)$ and the actions of $t$ and $\psi$ are twisted
by $\Chi$: 
\vglue12pt
\hfill $
\displaystyle{t(a) = -\omega^2 a = \omega^6 a \qquad\hbox{and}\qquad
\psi(a) = \omega a\ . }%
 $
\Endproof\vskip12pt 

We next write down the invariants
$E^{F}_\ast$ for the various finite subgroups $F$ of $G_{24}$. 
To do this, we work up from the symmetric algebra $S(\rho)$, and we use
the presentation of the
symmetric algebra as given in the exact sequence (3.9). 
Recall that we
have written $S(F) = \WW[x_1,x_2,x_3]$ where the normal subgroup of
order three in $G_{24}$ cyclically permutes the $x_i$. This action
by the cyclic group extends in an obvious way to an action of the
symmetric group $\Sigma_3$ on three letters; thus we have an inclusion
of algebras
$$
\WW[\sigma_1,\sigma_2,\sigma_3] = \WW[x_1,x_2,x_3]^{\Sigma_3}
\subseteq S(F)^{C_3}.
$$
There is at least one other obvious element invariant under the action
of $C_3$: set
\begin{equation}\label{eq100.6}
\epsilon = x_1^2x_2 + x_2^2x_3 + x_3^2x_1 - x_2^2x_1 - x_1^2x_3 - x_3^2x_2\ .
\end{equation}
This might be called the ``anti-symmetrization'' (with respect to $\Sigma_3$) 
of $x_1^2x_2$.

\begin{lem}\label{lem100.5} There is an isomorphism
$$
\WW[\sigma_1,\sigma_2,\sigma_3,\epsilon]/(\epsilon^2 - f) \cong
S(F)^{C_3}
$$
where $f$ is determined by the relation
$$
\epsilon^2 = -27\sigma_3^2 - 4\sigma_2^3 - 4\sigma_3\sigma_1^3 + 18
\sigma_1\sigma_2\sigma_3 + \sigma_1^2\sigma_2^2 \ .
$$
Furthermore{\rm ,} the actions of $t$ and $\psi$ are given by
$$
t(\sigma_1) = \omega^2\sigma_1\qquad t(\sigma_2) = -\sigma_2\qquad
t(\sigma_3) = \omega^6\sigma_3\qquad t(\epsilon) = \omega^2\epsilon
$$ 
and
$$
\psi(\sigma_1) = \omega\sigma_1\qquad \psi(\sigma_2) = \omega^2\sigma_2\qquad
\psi(\sigma_3) = \omega^3\sigma_3\qquad \psi(\epsilon) = \omega^3\epsilon\ .
$$
\end{lem}

\Proof  Except for the action of $\psi$, this is Lemma 4 of
\cite{GHM}. The action of $\psi$ is straightforward 
\Endproof\vskip4pt 

This immediately leads to the following result.

\begin{prop}\label{prop100.6} There is an isomorphism
$$
\WW[\sigma_2,\sigma_3,\epsilon]/(\epsilon^2 - g) \cong S(\rho)^{C_3}
$$
where $g$ is determined by the relation
$$
\epsilon^2 = -27\sigma_3^2 - 4\sigma_2^3
$$
with the actions of $t$ and $\psi$ as given above in Lemma {\rm \ref{lem100.5}.}
Under this isomorphism $\sigma_3$ maps to $d$.
\end{prop}

\Proof  This follows immediately from Lemma \ref{lem100.5},
the short exact sequence (\ref{eq100.5}), 
and the fact (see the proof of Lemma \ref{lem100.4}) 
that $H^1(C_3,S(F)) = 0$. Together these imply that
\vglue12pt
\hfill $
\displaystyle{S(\rho)^{C_3} \cong S(F)^{C_3}/(\sigma_1)\ .}
$ 
\Endproof\vskip12pt 

The next step is to invert the element $N$ of (\ref{eq100.4}).
This element is the image of $\sigma_3^4$; thus, we are effectively inverting
the element $d= \sigma_3 \in S(\rho)^{C_3}$. We begin with the observation
that if we divide
$$
\epsilon^2 = -27\sigma_3^2 - 4\sigma_2^3
$$
by $\sigma_3^6$ we obtain the relation
$$
(\frac{\epsilon}{\sigma_3^3})^2 + 4(\frac{\sigma_2}{\sigma_3^2})^3 =
-\frac{27}{\sigma^4_3}\ .
$$
Thus if we set 
\begin{equation}\label{eq100.7}
c_4 = - \frac{\omega^2\sigma_2}{\sigma_3^2},
\quad c_6 = \frac{\omega^3\epsilon}{2\sigma^3_3},
\quad \Delta = -\frac{\omega^6}{4\sigma_3^4} = \frac{\omega^2}{4\sigma_3^4}
\end{equation}
then we get the expected relation\footnote{This is the relation appearing
in theory of modular forms \cite{Deligne}, except here $2$ is invertible
so we can replace 1728 by 27. There is some discussion of the 
connection in \cite{GS}. 
The relation could be arrived at more naturally 
by choosing, as our basic formal group law, one arising from the 
theory of elliptic
curves, rather than the Honda formal group law.}
$$
c_6^2 - c_4^3 = 27\Delta\ .
$$
Furthermore, $c_4$, $c_6$, and $\Delta$ are all invariant under the action
of the entire group $G_{24}$. (Indeed, the powers of $\omega$ are introduced
so that this happens.)

To describe the group cohomology, we define elements
$$
\alpha = \frac{\omega a}{d} \in H^1(C_3,(S(\rho)[N^{-1}])_4)
$$
and
$$
\beta = \frac{\omega^3b}{d^2} \in H^2(C_3,(S(\rho)[N^{-1}])_{12}).
$$
These elements are fixed by $t$ and $\psi$ and, for degree reasons,
acted on trivially by $c_4$ and $c_6$. The following is now easy.

\begin{prop}\label{G24coho} {\rm 1 )} The inclusion
$$
\ZZ_3[c_4,c_6,\Delta^{\pm 1}]/(c_6^2 - c_4^3 = 27\Delta) \to
S(\rho)[N^{-1}]^{G_{24}}
$$
is an isomorphism of algebras.
\smallbreak 
{\rm 2)} There is an exact sequence
$$
S(\rho)[N^{-1}] \mathop{\longrightarrow}^{\tran} H^\ast(G_{24},
S(\rho)[N^{-1}]) \to \FF_3[\alpha,\beta,\Delta^{\pm 1}]/(\alpha^2) \to 0
$$
and $c_4$ and $c_6$ act trivially on $\alpha$ and $\beta$.
\end{prop}

Then a completion argument, as in Theorem 6 of \cite{GHM} or
\cite{Nave} implies the next result.

\begin{thm}\label{G24cohocomplete} {\rm 1)} There is an isomorphism of algebras
$$
(E_\ast)^{G_{24}} \cong
\ZZ_3[[c_4^3\Delta^{-1}]][c_4,c_6,\Delta^{\pm 1}]/(c_6^2 - c_4^3 = 27\Delta)\ .
$$

{\rm 2)} There is an exact sequence
$$
E_\ast \mathop{\longrightarrow}^{\tran} H^\ast(G_{24},
E_\ast) \to \FF_3[\alpha,\beta,\Delta^{\pm 1}]/(\alpha^2) \to 0
$$
and $c_4$ and $c_6$ act trivially on $\alpha$ and $\beta$.
\end{thm}

\begin{rem}\label{othergroups}
The same kind of reasoning can be used to obtain the group
cohomologies $H^\ast(F,E_\ast)$ for other finite subgroups
of $\GG_{2}$. First define an element
$$
\delta = \sigma_3^{-1} \in S(\rho)[N^{-1}]\ .
$$
Then $\Delta = (\omega^2/4)\delta^4$; thus $-\Delta$ has a
square root:
$$
(-\Delta)^{1/2} = \frac{\omega^3}{2}\delta^2\ .
$$
The elements $t$ and $\psi$ of $G_{24}$ act on $\delta$ by
the formulas
$$
t(\delta) = \omega^2\delta\qquad\hbox{and}\qquad \psi(\delta) =
\omega^5\delta\ .
$$
The element $(-\Delta)^{1/2}$ is invariant under the action of
$t^2$ and $\psi$ (whereas the evident square root of $\Delta$
is not fixed by $\psi$).

Let $C_{12}$ be the cyclic subgroup of order $12$ in $G_{24}$ generated
by $s$ and $\psi$.
This subgroup has a cyclic subgroup $C_6$ of order $6$ generated
by $s$ and $t^2 = \psi^2$. We have
\begin{align*}
(E_\ast)^{C_3}& \cong \WW[[c_4^3\Delta^{-1}]][c_4,c_6,\delta^{\pm 1}]/
(c_6^2 - c_4^3 = 27\Delta)\\
(E_\ast)^{C_{12}}& \cong
\ZZ_3[[c_4^3\Delta^{-1}]][c_4,c_6,(-\Delta)^{\pm 1/2}]/(c_6^2
- c_4^3 = 27\Delta)\\
(E_\ast)^{C_6} &\cong \WW \otimes_{\ZZ_3} (E_\ast)^{C_{12}}\\
(E_\ast)^{G_{12}}& \cong \WW[[c_4^3\Delta^{-1}]][c_4,c_6,\Delta^{\pm 1}]/
(c_6^2 - c_4^3 = 27\Delta)\cong \WW \otimes_{\ZZ_3} (E_\ast)^{G_{24}}.
\end{align*}
Furthermore, for all these groups, the analogue of Theorem
\ref{G24cohocomplete}.2 holds. For example, there are exact sequences
\begin{align*}
E_\ast &\mathop{\longrightarrow}^{\tran} H^\ast(C_3,E_\ast) 
\to \FF_9[\alpha,\beta,\delta^{\pm 1}]/(\alpha^2) \to 0\\
E_\ast &\mathop{\longrightarrow}^{\tran} H^\ast(C_{12},E_\ast)
\to \FF_3[\alpha,\beta,(-\Delta)^{\pm 1/2}]/(\alpha^2) \to 0\\
E_\ast &\mathop{\longrightarrow}^{\tran} H^\ast(C_{6},E_\ast)
\to \FF_9[\alpha,\beta,(-\Delta)^{\pm 1/2}]/(\alpha^2) \to 0\\
E_\ast &\mathop{\longrightarrow}^{\tran} H^\ast(G_{12},E_\ast)
\to \FF_9[\alpha,\beta,(\Delta)^{\pm 1}]/(\alpha^2) \to 0 
\end{align*}
and $c_4$ and $c_6$ act trivially on $\alpha$ and $\beta$.
\end{rem}

These results allow one to competely write down the various homotopy fixed
point spectral sequences for computing $\pi_\ast E^{hF}$ for
the various finite groups in question. The differentials 
in the spectral sequence follow from
Toda's classical results and the following easy
observation: every element in the image of the transfer is a permanent
cycle. We record:

\begin{lem}\label{lem100.12} In the spectral sequence
$$
H^\ast (G_{24},E_\ast) \Longrightarrow \pi_\ast E^{hG_{24}}
$$
the only nontrivial differentials are $d_5$ and $d_9$.  They  
are determined by 
$$
d_5(\Delta) = a_1\alpha\beta^2\qquad\hbox{and}\qquad 
d_9(\alpha \Delta^2) = a_2 \beta^5 
$$
where $a_1$ and $a_2$ are units in $\FF_3$. 
\end{lem}

\Proof  These are a consequence of Toda's famous differential
(see \cite{Toda}) and nilpotence. See Proposition 7 of \cite{GHM} or,
again, \cite{Nave}. There it is done for $G_{12}$ rather than $G_{24}$,
but because $G_{12}$ is of index $2$ in $G_{24}$ and 
we are working at the prime $3$, this is sufficient.
\Endproof\vskip4pt

The lemma immediately calculates the 
differentials in the other spectral
sequences; for example, if one wants homotopy fixed points with respect
to the $C_3$-action, we have, up to units,
$$
d_5(\delta) = \delta^{-3}\alpha\beta^2\qquad\hbox{and}\qquad
d_9(\alpha\delta^2) = \delta^{-6}\beta^5.
$$

It is also worth pointing out that the $d_5$-differential in Lemma 
\ref{lem100.12} and some standard Toda bracket manipulation (see 
the proof of Theorem 8 in \cite{GHM}) implies the relation 
$(\Delta\alpha)\alpha= \pm \beta^3$ which holds in $\pi_{27}(E^{hG_{24}})$. 

The above discussion is summarized in the following main homotopy 
theoretic result of this section.

\begin{thm}\label{thm100.13} In the spectral sequence
$$
H^\ast (C_3,E_\ast) \Longrightarrow \pi_\ast E^{hC_3}
$$
we have an inclusion of subrings
$$
E_\infty^{0,\ast} \cong \WW[[c_4^3\Delta^{-1}]][c_4,c_6,c_4\delta^{\pm 1}, 
c_6\delta^{\pm 1}, 3\delta^{\pm 1},\delta^{\pm 3}]/(c_4^3-c_6^2 = 27\Delta) 
\subseteq E_2^{0,\ast}.
$$
In positive filtration $E_\infty$ is additively generated by the
elements $\alpha${\rm ,} $\delta\alpha${\rm ,} $\alpha\beta${\rm ,} $\delta\alpha\beta${\rm ,}
$\beta^j${\rm ,} $1 \leq j \leq 4$ and all multiples of these elements
by $\delta^{\pm 3}$. These elements are of order $3$
and are annihilated by $c_4${\rm ,} $c_6${\rm ,} $c_4\delta^{\pm 1}${\rm ,} $c_6\delta^{\pm 1}$ 
and $3\delta^{\pm 1}$. Furthermore the following multiplicative 
relation holds in $\pi_{30}(E^{hC_3})${\rm :} 
$\delta^3(\delta\alpha)\alpha=\pm\omega^6\beta^3$. 
\end{thm}

For the case of the cyclic group $C_6$ of order $6$ generated
by $s$ and $t^2$, one notes that $t^2(\delta) = -\delta$ and the
spectral sequence can now be read off Theorem~\ref{thm100.13}. 
This also determines the case of $C_{12}$, the group 
generated by $s$ and~$\psi$. We leave 
the details to the reader but state the result in the case of~$G_{24}$. 

\begin{thm}\label{thm100.14} In the spectral sequence
$$
H^\ast (G_{24},E_\ast) \Longrightarrow \pi_\ast E^{hG_{24}}
$$
we have an inclusion of subrings
$$
E_\infty^{0,\ast} \cong \ZZ_3[[c_4^3\Delta^{-1}]][c_4,c_6,
c_4\Delta^{\pm 1}, c_6\Delta^{\pm 1}, 3\Delta^{\pm 1}, 
\Delta^{\pm 3}]/(c_4^3-c_6^2 = 27\Delta) \subseteq E_2^{0,\ast}.
$$
In positive filtration $E_\infty$ is additively generated by the
elements $\alpha${\rm ,} $\Delta\alpha${\rm ,} $\alpha\beta${\rm ,} $\Delta\alpha\beta${\rm ,}
$\beta^j${\rm ,} $1 \leq j \leq 4$ and all multiples of these elements
by $\Delta^{\pm 3}$. These elements are of order $3$
and are annihilated by $c_4${\rm ,} $c_6${\rm ,} $c_4\Delta^{\pm 1}${\rm ,} $c_6\Delta^{\pm 1}$ 
and $3\Delta^{\pm 1}$. Furthermore $(\Delta\alpha)\alpha=\pm \beta^3$
in $\pi_{30}(E^{hG_{24}})$.
\end{thm}

\begin{rem}\label{periodicity} {\rm 1)} Note that $E^{hC_3}$ is periodic of
period $18$ and that
$\delta^3$ detects a periodicity class. The spectra $E^{hC_6}$
and $E^{hC_{12}}$ are periodic with period $36$ and
$(-\Delta)^{3/2}$ detects the periodicity generator. Finally
$E^{hG_{12}}$ and $E^{hG_{24}}$ are periodic with period
$72$ and $\Delta^3$ detects the periodicity generator.
\smallbreak 
{\rm 2)} By contrast, the Morava module $E_\ast E^{hG_{24}}$ 
is of period $24$. 
To see this, note that the isomorphism of (\ref{fixediso}) 
supplies an isomorphism of Morava modules
$$
E_\ast E^{hG_{24}} \cong \Hom^c(\GG_2/G_{24},E_\ast) 
$$
with $\GG_2$ acting diagonally on the right-hand side.  
Then $G_{24}$-invariance of $\Delta$ implies that there is a well 
defined automorphism of Morava modules given by 
$\varphi\mapsto (g\mapsto \varphi(g)g_*(\Delta))$.   

3) If $F$ has order prime to $3$, then 
$\pi_\ast(E^{hF}) = (E_\ast)^{F}$ is easy to calculate. For example,
using (\ref{action-of-fstar}) and (\ref{devhop1}) 
we obtain
$$
\pi_\ast (E^{hSD_{16}}) = \ZZ_3[v_1][v_2^{\pm 1}]{\ }\widehat{ }\ 
\subseteq E_\ast
$$
with $v_1=u_1u^{-2}$ and $v_2=u^{-8}$ and completion is with 
respect to the ideal generated by $(v_1^4v_2^{-1})$. Similarly 
$$
\pi_\ast (E^{hQ_8}) = \ZZ_3[v_1][\omega^2u^{\pm 4}]{\ }\widehat{ }\ 
\subseteq E_\ast 
$$
where completion is with respect to the ideal $(v_1^2\omega^2u^4)$. 
Note that $E^{hSD_{16}}$ is periodic of order $16$ and 
$E^{hQ_8}$ is periodic of order $8$. 
\end{rem}

We finish this section by listing exactly the computational results 
we will use in building the tower in Section 5.

\begin{cor}\label{iso-in-degree-zero} {\rm 1)} Let $F \subseteq \GG_2$
be any finite subgroup. Then the edge homomorphism
$$
\pi_0E^{hF} \longrightarrow (E_0)^F
$$
is an isomorphism of algebras.

{\rm 2)} Let $F \subseteq \GG_2$ be any finite subgroup. Then
$$
\pi_{24}E^{hF} \longrightarrow (E_{24})^F
$$
is an injection.
\end{cor}

\Proof  If $F$ has order prime to $3$, both of these statements
are clear. If $3$ divides the order of $F$, then the 3-Sylow subgoup
of $F$ is conjugate to $C_3$; hence $\pi_\ast E^{hF}$ is a retract of
$\pi_\ast E^{hC_3}$ and the result follows from
Theorem \ref{thm100.13}. \phantom{overthere}
\hfill\qed

\begin{cor}\label{zero-in-pi1} Let $F \subseteq \GG_2$ be a finite
subgroup. If the order of $F$ is prime to $3${\rm ,} then
$$
\pi_1E^{hF} = 0\ .
$$
For all finite subgroups $F${\rm ,} 
$$
\pi_{25}E^{hF} = 0\ .
$$
\end{cor}

\Proof  Again apply Theorem \ref{thm100.13}.
\hfill\qed

\begin{cor}\label{effect-of-the-center} Let $F \subseteq \GG_2$ be
any finite subgroup containing the central element $\omega^4 = -1$.
Then 
$$
\pi_{26}E^{hF} = 0\ .
$$
\end{cor}

\Proof  Equation \ref{action-of-fstar} implies that
$(E_\ast)^F$ will be concentrated in degrees congruent to $0$  
mod $4$. Combining this observation with Theorem \ref{thm100.13}
proves the result.
\hfill\qed

\section{The algebraic resolution}

Let $G$ be a profinite group. Then its $p$-adic group ring $\ZZ_p[[G]]$ 
is defined as $\lim_{U,n}\ZZ_p/(p^n)[G/U]$ 
where $U$ runs through all open subgroups of $G$. Then $\ZZ_p[[G]]$ is a 
complete ring and we will only consider continuous modules over such rings.

In this section we will construct 
our resolution of the trivial $\ZZ_3[[\GG_2]]$-module $\ZZ_3$. 
In \ref{rem0.3} we wrote down a splitting of the 
group $\GG_2$ as $\GG_2^1 \times
\ZZ_3$ and this splitting allows us to focus on constructing 
a resolution of $\ZZ_3$ as a $\ZZ_3[[\GG_2^1]]$-module.

Recall that we have selected a maximal finite subgroup $G_{24}
\subseteq \GG_2^1$; it is generated by an element $s$ of order three,
$t=\omega^2$, and $\psi = \omega\phi$ where $\omega$ is a primitive eighth
root of unity and $\phi$ is the Frobenius. As before $C_3$ denotes the normal
subgroup of order $3$ in $G_{24}$ generated by $s$, 
$Q_8$ denotes the subgroup of $G_{24}$ of order $8$
generated by $t$ and $\psi$ and $SD_{16}$ the subgroup of $\GG_2^1$ 
generated by $\omega$ and $\psi$. 
\def\twisted{{{\chi\uparrow^{\GG_2^1}_{SD_{16}}}}}
\def\largetwisted{{{\chi\uparrow^{\GG_2}_{SD_{16}}}}}

The group $Q_8$ is a subgroup of $SD_{16}$ of index $2$.
Let $\chi$ be the sign representation (over $\ZZ_3$)
of $SD_{16}/Q_8$; we regard $\chi$ as representation of $SD_{16}$ using
the quotient map. (Note that this is not the same $\chi$ as in \S 3!)  
In this section, the induced $\ZZ_3[[\GG_2^1]]$-module
\begin{equation}\label{twisted}
\twisted \defeq \ZZ_3[[\GG^1_2]] \otimes_{\ZZ_3[SD_{16}]} \chi
\end{equation}
will play an important role. If $\tau$ is the trivial representation of 
$SD_{16}$,
there is an isomorphism of $SD_{16}/Q_8$-modules
\begin{equation}\label{twisted-plus-trivial}
\ZZ_3[SD_{16}/Q_8] \cong \chi \oplus \tau\ .
\end{equation}
Thus, we have that $\twisted$ is a direct summand of the 
induced module 
$$\ZZ_3[[\GG_2^1/Q_8]]=\ZZ_3[[\GG_2^1]]\otimes_{\ZZ_3[Q_8]}\ZZ_3 \ .  
$$

The following is the main algebraic result of the paper. It will
require the entire section to prove.

\begin{thm}\label{thm4.1} There is an exact sequence of 
$\ZZ_3[[\GG_2^1]]$\/{\rm -}\/modules
$$
0 \to \ZZ_3[[\GG_2^1/G_{24}]] \to \twisted \to
\twisted \to \ZZ_3[[\GG_2^1/G_{24}]] \to \ZZ_3 \to 0\ .
$$
\end{thm}
\vskip9pt

Some salient features of this ``resolution'' (we will use this word
even though not all of the modules are projective) 
are that each module is a summand of a permutation module $\ZZ_3[[\GG_2^1/F]]$ 
for a finite subgroup $F$ and that each module is free over $K$, 
where $K \subseteq \GG^1_2$ is a subgroup so that
we can decompose the 3-Sylow subgroup $S_2^1$ of $\GG_2$ as $K\rtimes C_3$.
Important features of $K$ include that it is a torsion-free $3$-adic
Poincar\'e duality group of dimension $3$. (See before Lemma \ref{lem4.10}
for more on $K$.)

Since $\GG_2 \cong \GG_2^1 \times \ZZ_3$, we may tensor the resolution
of Theorem \ref{thm4.1} with the standard resolution for $\ZZ_3$ to get
the following as an immediate corollary,

\begin{cor}\label{cor4.2} There is an exact sequence of 
$\ZZ_3[[\GG_2]]$\/{\rm -}\/modules
\begin{align*}
0 \to \ZZ_3[[\GG_2/G_{24}]] &\to \ZZ_3[[\GG_2/G_{24}]] \oplus 
\largetwisted \to \largetwisted \oplus \largetwisted\\
& \to \largetwisted \oplus
\ZZ_3[[\GG_2/G_{24}]] \to \ZZ_3[[\GG_2/G_{24}]] \to \ZZ_3 \to 0\ .
\end{align*}
\end{cor}
\vskip8pt

As input for our calculation we will use 
$$H^\ast(S_2^1):=H^\ast(S_2^1,\FF_3)=\Ext^\ast_{\ZZ_3[[S_2^1]]}(\ZZ_3,\FF_3)\ , 
$$
as calculated by the second author in
\cite{Henn}. This is an effective starting point
because of the following lemma and the fact that
\begin{equation}\label{dual}
\Ext_{\ZZ_3[[S^1_2]]}^q(M,\FF_3) \cong
\Tor^{\ZZ_3[[S^1_2]]}_q(\FF_3,M)^\ast 
\end{equation} 
for any profinite continuous $\ZZ_3[[S_2^1]]$-module $M$. 
Here $(-)^\ast$ means $\FF_p$-linear dual.

A profinite group $G$ is called finitely generated if there is
a finite set of elements $X \subseteq G$ so that the subgroup
generated by $X$ is dense. This is true of all the groups in
this paper. If $G$ is a $p$-profinite group and $I \subseteq \ZZ_p[[G]]$
is the kernel of the augmentation $\ZZ_p[[G]] \to \FF_p$, then
$$
\ZZ_p[[G]]  \cong \lim_n \ZZ_p[[G]]/I^n.
$$
A $\ZZ_p[[G]]$-module $M$ will be called {\it complete} if it is $I$-adically
complete, i.e.\ if $M\cong \lim_nM/I^nM$. 

\begin{lem}\label{lem4.3} \hskip-4pt Let $G$ be a finitely generated
$p$\/{\rm -}\/profinite group and \hbox{$f:M \to N$} a morphism of complete 
$\ZZ_p[[G]]$\/{\rm -}\/modules. If
$$
\FF_p \otimes f: \FF_p \otimes_{\ZZ_p[[G]]} M \to \FF_p \otimes_{\ZZ_p[[G]]} N
$$
is surjective{\rm ,} then $f$ is surjective. If
$$
\Tor(\FF_p,f): \Tor^{\ZZ_p[[G]]}_q(\FF_p,M) \to \Tor^{\ZZ_p[[G]]}_q(\FF_p,N)
$$
is an isomorphism for $q = 0$ and onto for $q=1${\rm ,} then $f$ is an
isomorphism.
\end{lem}

\Proof  This is an avatar of Nakayama's lemma.
To see this, suppose $K$ is some complete $\ZZ_p[[G]]$-module so that
$\FF_p \otimes_{\ZZ_p[[G]]} K = 0$. Then an inductive argument shows
$$
\ZZ_p[[G]]/I^n \otimes_{\ZZ_p[[G]]} K = 0
$$
for all $n$; hence $K = 0$. This is the form of Nakayama's lemma
we need. 

The result is then proved using the long exact sequence of $\Tor$ groups:
the weaker hypothesis implies that the cokernel of $f$ is trivial; the
stronger hypothesis then implies that the kernel of $f$ is trivial.
\Endproof\vskip4pt 

We next turn to the details about $H^\ast(S_2^1;\FF_3)$ from \cite{Henn}. (See
Theorem 4.3 of that paper.) We will omit the coefficients $\FF_3$ in order 
to simplify our notation. 
The key point here is that the cohomology
of the group $S^1_2$ is detected on the centralizers of the cyclic
subgroups of order $3$. There are two conjugacy classes of
such subgroups of order $3$ in $S_2^1$; namely, $C_3$ and $\omega C_3
\omega^{-1}$. The element $s = s_1$ is our chosen generator for $C_3$;
thus we choose as our generator for $\omega C_3\omega^{-1}$ the
element $s_2 = \omega s \omega^{-1}$. The Frobenius $\phi$ also
conjugates $C_3$ to $\omega C_3 \omega^{-1}$ and a short
calculation shows that
\begin{equation}\label{frobaction1}
\phi (s_1) = \phi (s)  = s_2^2\ .
\end{equation}
The centralizer $C(C_3)$ in $S^1_2$ is isomorphic to $C_3 \times \ZZ_3$,  
$\omega\phi$ commutes with $C(C_3)$ and conjugation by 
$\omega^2$ sends $x\in C(C_3)$ to its inverse $x^{-1}$ (see \cite{GHM}). 
In particular, for every $x \in C(C_3)$ we have
\begin{equation}\label{frobaction2}
\omega x \omega^{-1} = \phi( x)^{-1}\phi^{-1} \in C(\omega C_3\omega^{-1})\ .
\end{equation}
Note that $C(\omega C_3 \omega^{-1}) = \omega C(C_3) \omega^{-1}$.
Write $E(X)$ for the exterior algebra on a set $X$. Then
$$
H^\ast(C(C_3)) \cong \FF_3[y_1] \otimes E(x_1,a_1)
$$
and
$$
H^\ast(C(\omega C_3\omega^{-1}))\cong  \FF_3[y_2] \otimes E(x_2,a_2)\ .
$$
We know that $C_4$ (which is generated by $\omega^2$) acts on 
$$
H^\ast(C(C_3))\cong \FF_3[y_1] \otimes E(x_1,a_1)
$$
sending all three generators to their negative. This action extends
to an action of $SD_{16}$ on the product 
$$
H^\ast(C(C_3)) \times H^\ast(C(\omega C_3\omega^{-1}))
$$
as follows. By (\ref{frobaction1}) and (\ref{frobaction2}) the action 
of the generators $\omega$ and $\phi$ of $SD_{16}$ 
is given by
\begin{alignat}{6}\label{cohoaction}
& \omega_\ast(x_1) = x_2,& \quad \omega_\ast(y_1) &= y_2,&\quad \omega_\ast(a_1)&=
a_2,\\
&\phi_\ast(x_1) = -x_2,& \quad \phi_\ast(y_1) &= -y_2,&\quad \phi_\ast(a_1)
&= -a_2.\nonumber
\\[6pt]
 \label{cohoaction1}
&\omega_\ast(x_2) = -x_1,&\quad \omega_\ast(y_2)& = -y_1, 
&\quad \omega_\ast(a_2)&=-a_1,\\
& \phi_\ast(x_2) = -x_1,& \quad \phi_\ast(y_2) &= -y_1,&\quad \phi_\ast(a_2)
&= -a_1\ .\nonumber
\end{alignat}
\vglue8pt

\begin{thm}[\cite{Henn}]\label{thm4.4} {\rm 1)} The inclusions
$\omega^i C(C_3)\omega^{-i} \to S_2^1${\rm ,} $i = 0,1$ induce an 
$SD_{16}$\/{\rm -}\/equivariant homomorphism
$$
H^\ast(S_2^1) \to \prod_{i=1}^2\FF_3[y_i] \otimes E(x_i,a_i) 
$$
which is an injection onto the subalgebra generated by
$x_1$\/{\rm ,}\/ $x_2${\rm ,} $y_1${\rm ,} $y_2${\rm ,} $x_1a_1 - x_2a_2${\rm ,} $y_1a_1$ 
and $y_2a_2$.
\smallskip

{\rm 2)} In particular{\rm ,} $H^\ast(S_2^1)$ is free as a module over
$\FF_3[y_1 + y_2]$ on generators $1${\rm ,} $x_1${\rm ,} $x_2${\rm ,} $y_1${\rm ,} 
$x_1a_1 - x_2a_2${\rm ,} $y_1a_1${\rm ,} $y_2a_2${\rm ,} and $y_1x_1a_1$.
\end{thm}

We will produce the resolution of Theorem \ref{thm4.1} from this
data and by splicing together the short exact sequences
of Lemma \ref{lem4.5}, \ref{lem4.6}, and \ref{lem4.7} below. Most
of the work will be spent in identifying the last module; this
is done in Theorem \ref{thm4.9}.

\def\bigcoh{{\Ext(\ZZ_3) \odot \FF_3[SD_{16}]}}

In the following computations, we will write
$$
\Ext(M) = \Ext_{\ZZ_3[[S^1_2]]}^\ast(M,\FF_3)\ .
$$
This graded vector space is a module over
$$
H^\ast(S^1_2) = \Ext_{\ZZ_3[[S^1_2]]}^\ast(\ZZ_3,\FF_3) = \Ext(\ZZ_3)\ ,
$$
and, hence is also a module over the sub-polynomial algebra
of $H^\ast(S^1_2)$ generated by $y_1 + y_2$. If $M$ is actually
a continuous $\ZZ_3[[\GG_2^1]]$-module, then $\Ext(M)$ has an action
by $SD_{16}\cong \GG_2^1/S_2^1$ which extends the action by $H^\ast(S^1_2)$
in the obvious way: if $\alpha\in SD_{16}$, 
$a \in H^\ast(S^1_2)$ and $x \in \Ext(M)$, then
$$
\alpha(ax) = \alpha(a)\alpha(x)\ .
$$
We can write this another way. Let's define $\bigcoh$ to be
the algebra constructed by taking
$$
\Ext(\ZZ_3) \otimes \FF_3[SD_{16}]
$$
with twisted product
$$
(a \otimes \alpha)(b \otimes \beta) = a\alpha(b) \otimes \alpha\beta\ .
$$
The above remarks imply that if $M$ is a $\ZZ_3[[\GG^1_2]]$-module, then 
$\Ext(M)$ is an $\bigcoh$-module.

This structure behaves well with respect to long exact sequences. If
$$
0 \to M_1 \to M_2 \to M_3 \to 0
$$
is a short exact sequence of $\ZZ_3[[\GG^1_2]]$-modules we get
a long exact sequence in $\Ext$ which is a long exact
sequence of $\bigcoh$-modules. As a matter of notation, if $x\in\Ext(\ZZ_3)$
we will write $\overline{x}\in\Ext(M)$ if $x$ is the
image of $\overline{x}$ under some unambiguous and injective sequence
of boundary homomorphisms of long exact sequences.


\begin{lem}\label{lem4.5} There is a short exact sequence of
$\ZZ_3[[\GG_2^1]]$-modules
$$
0 \to N_1 \to \ZZ_3[[\GG_2^1/G_{24}]] \mathop{\longrightarrow}^{\epsilon}
\ZZ_3\to 0
$$
where the map $\epsilon$ is the augmentation. If we write $z$ 
for $x_1a_1 - x_2a_2 \in H^2(S_1^2)$ then 
$\Ext(N_1)$ is a module over $\FF_3[y_1+y_2]$ on generators $e${\rm ,}
$\overline{z}${\rm ,} $\overline{y_1a_1}${\rm ,} $\overline{y_2a_2}${\rm ,}
and $\overline{y_1x_1a_1}$
of degrees $0${\rm ,} $1${\rm ,} $2${\rm ,} $2${\rm ,} and $3$ respectively. The last four
generators are free and $(y_1 + y_2)e = 0$. The action of $SD_{16}$
is determined by the action on $\Ext(\ZZ_3)$ and the facts that
$$
\omega_\ast(e) = -e = \phi_\ast(e)\ .
$$
\end{lem}

\Proof  As a $\ZZ_3[[S^1_2]]$-module, there is an
isomorphism
$$
\ZZ_3[[\GG_2^1/G_{24}]] \cong \ZZ_3[[S^1_2/C_3]] \oplus
\ZZ_3[[S^1_2/\omega C_3\omega^{-1}]].
$$
Hence, by the Shapiro lemma there is an isomorphism
$$
\Ext(\ZZ_3[[\GG_2^1/G_{24}]]) \cong H^\ast(C_3,\FF_3) \times
H^\ast(\omega C_3\omega^{-1},\FF_3)
$$
and the map $\Ext(\ZZ_3) \to \Ext(\ZZ_3[[\GG_2^1/G_{24}]])$
corresponds via this isomorphism to the restriction map. 
The result now follows from Theorem \ref{thm4.4}. 
\Endproof\vskip4pt 

Recall that $\chi$ is the rank one (over $\ZZ_3$) representation of
$SD_{16}$ obtained by pulling back the sign representation along
the quotient map $\varepsilon:SD_{16} \to SD_{16}/Q_8 \cong \ZZ/2$.

\begin{lem}\label{lem4.6} There is a short exact sequence of
$\ZZ_3[[\GG_2^1]]$-modules
$$
0 \to N_2 \to \twisted \to N_1 \to 0\ .
$$
The cohomology module $\Ext(N_2)$ is a freely generated
module over $\FF_3[y_1+y_2]$ on generators $\overline{z}${\rm ,} $\overline{y_1a_1}${\rm ,}
$\overline{y_2a_2}${\rm ,} and $\overline{y_1x_1a_1}$
of degrees $0${\rm ,} $1${\rm ,} $1${\rm ,} and $2$ respectively. 
The action of $\omega$ is determined by the action on $\Ext(\ZZ_3)$.
\end{lem}

\Proof  By the previous result the $SD_{16}$-module
$\FF_3 \otimes_{\ZZ_3[[S^1_2]]} N_1$ is one dimensional over
$\FF_3$ generated by the dual (with respect to (\ref{dual})) 
of the class of $e$ and the action is given 
by the sign representation along $\varepsilon$.
Lift $e$ to an element $d \in N_1$. Then $SD_{16}$ may not
act correctly on $d$, but we can average $d$ to obtain an element
$c$ on which $SD_{16}$ acts correctly and which reduces to the same
element in $\FF_3 \otimes_{\ZZ_3[[S^1_2]]} N_1$; indeed,
$$
c = \frac{1}{16}\sum_{\alpha \in SD_{16}} 
\varepsilon(\alpha)^{-1}\alpha_\ast (d)\ .
$$
This defines the morphism
$$
\twisted \to N_1\ .
$$
Lemma \ref{lem4.3} now implies that this map is surjective and we obtain
the exact sequence \pagebreak we need. For the calculation of $\Ext(N_2)$ note that
we have an isomorphism of $S^1_2$-modules
$$
\twisted \cong \ZZ_3[[S^1_2]]\ .
$$
The result now follows from the previous lemma and the long exact sequence.
\phantom{overthere}\hfill\qed

\begin{lem}\label{lem4.7} There is a short exact sequence of
$\ZZ_3[[\GG_2^1]]$\/{\rm -}\/modules
$$
0 \to N_3 \to \twisted \to N_2 \to 0
$$
where $\Ext(N_3)$ is a free module over $\FF_3[y_1+y_2]$ on generators
$\overline{y_1a_1}${\rm ,} $\overline{y_2a_2}${\rm ,} $\overline{y_1x_1a_1}${\rm ,}
and $\overline{y_2x_2a_2}$ of degree $0${\rm ,} $0${\rm ,} $1$ and $1$ respectively. 
In fact{\rm ,} the iterated boundary homomorphisms
$$
\Ext^\ast(N_3) \to \Ext^{\ast+3}(\ZZ_3) = H^{\ast + 3}(S^1_2,\FF_3)
$$
define an injection onto an $\bigcoh$\/{\rm -}\/submodule isomorphic to
$$\Ext(\ZZ_3[[\GG^1_2/G_{24}]])\ .$$
\end{lem}

\Proof  The $SD_{16}$-module $\FF_3 \otimes_{\ZZ_3[[S^1_2]]} N_2$
is $\FF_3 \otimes_{\ZZ_3} \chi$ generated by the class dual 
to $\bar{z}$.  As in the
proof the last lemma, we can now form a surjective map
$$
\twisted \to N_2\ ,
$$
and this map defines $N_3$. The calculation of $\Ext(N_3)$ follows from
the long exact sequence.
\Endproof\vskip4pt 

To make use of this last result we prove a lemma.

\begin{lem}\label{lem4.8} Let $A = \bigcoh$ and $M =
\Ext(\ZZ_3[[\GG^1_2/G_{24}]])${\rm ,} regarded as an $A$-module. Then
$M$ is a simple $A$-module\/{\rm ;}\/ in fact{\rm ,}
$$
{\rm End}_A(M) \cong \FF_3\ .
$$
\end{lem}

\Proof  Let $e_1$ and $e_2$ in
$$
\Ext(\ZZ_3[[\GG^1_2/G_{24}]]) \cong H^\ast(C_3,\FF_3) \times H^\ast(\omega
C_3 \omega^{-1},\FF_3)
$$
be the evident two generators in degree $0$. If $f: M\to M$ is any $A$-module
endomorphism, we may write
$$
f(e_1) = ae_1 + be_2
$$
where $a,b \in \FF_3$. Then (using the notation of Theorem \ref{thm4.4})
we have
$$
0 = f(y_2e_1) = by_2e_2.
$$
Since $y_2e_2 \ne 0$, we have $b=0$. Also, since $\omega_\ast(e_1) = e_2$,
we have $f(e_2) = ae_2.$ Finally, since every homogeneous element of
$M$ is of the from $x_1y^i_1e_1 + x_2y_2^ie_2$, we have $f = a\hbox{id}_M$.
\Endproof\vskip4pt 

This means that in order to prove the following result, we need only
produce a map $f:N_3 \to \ZZ_3[[\GG_2^1/G_{24}]]$ of 
$\GG^1_2$-modules which induces a nonzero map on $\Ext$ groups.

\begin{thm}\label{thm4.9} There is an isomorphism of 
$\ZZ_3[[\GG_2^1]]$-modules
$$
N_3 \mathop{\longrightarrow}^{\cong} \ZZ_3[[\GG_2^1/G_{24}]]\ .
$$
\end{thm}

This requires a certain amount of preliminaries, and some further lemmas. 
We are looking for a diagram (see diagram (\ref{eq4.5}) below) which will
build and detect the desired map.

The first ingredient of our calculation is a spectral sequence.
Let us write
$$
0 \to C_3 \to C_2 \to C_1 \to C_0 \to \ZZ_3 \to 0
$$
for the resolution obtained by splicing together the short
exact sequences of Lemma \ref{lem4.5}, \ref{lem4.6}, and \ref{lem4.7}:
$$
0 \to N_3 \to \twisted \to \twisted \to
\ZZ_3[[\GG^1_2/G_{24}]] \to \ZZ_3 \to 0\ .
$$
By extending the resolution $C_\bullet \to \ZZ_3$ to a
bicomplex of projective $\ZZ_3[[\GG^1_2]]$-modules, we get, for
any $\ZZ_3[[\GG^1_2]]$-module $M$ and any closed subgroup 
$H \subseteq \GG^1_2$,
a first quadrant cohomology spectral sequence
\begin{equation}\label{eq4.1}
E^{p,q}_1 = \Ext_{\ZZ_3[[H]]}^p(C_q,M) \Longrightarrow H^{p+q}(H,M)\ .
\end{equation}
In particular, because $E_1^{0,q} = 0$ for $q > 3$,  there is 
an edge homomorphism
\begin{equation}\label{eq4.2}
\Hom_{\ZZ_3[[H]]}(N_3,M) = \Hom_{\ZZ_3[[H]]}(C_3,M) \to H^3(H,M)\ .
\end{equation}
Dually, there are homology spectral sequences
\begin{equation}\label{eq4.3}
E_{p,q}^1 = \Tor^{\ZZ_3[[H]]}_p(M,C_q) 
\Longrightarrow H_{p+q}(H,M)
\end{equation}
with an edge homomorphism
\begin{equation}\label{eq4.4}
H_3(H,M) \to M \otimes_{\ZZ_3[[H]]} N_3\ .
\end{equation}

That said, we remark that the important ingredient here
is that $\GG^1_2$ contains a subgroup $K$
which is a Poincar\'e duality group of dimension three and which
has good cohomological properties. 
The reader is referred to \cite{SW} for a modern discussion of 
a duality theory in the cohomology of profinite groups.

To define $K$, we use the filtration on the $3$-Sylow subgroup 
$S^1_2 = F_{1/2}\SS^1_2$ of $\GG_2^1$ described in the
first section. There is a projection
$$
S^1_2 \to  F_{1/2}\SS^1_2/ F_{1}\SS^1_2 \cong \FF_9\ .
$$
We follow this by the map $\FF_9\to \FF_9/\FF_3\cong C_3$ 
to define a group homorphism $S^1_2 \to C_3$; then, we define
$K \subset S^1_2$ to be the kernel. The chosen subgroup
$C_3 \subseteq S^1_2$ of order $3$ provides a splitting
of $S^1_2 \to C_3$; hence $S^1_2$ can be written as a semi-direct
product $K\rtimes C_3$. Note that every element of order three in $S^1_2$ maps
to a nonzero element in $C_3$ so that $K$ is torsion free. 

From \cite{Henn}, we know a good deal about $K$, some of which
is recorded in the following lemma. Let $j:K \to S^1_2$ denote the
inclusion. 

\begin{lem}\label{lem4.10} The group $K$ is a $3$\/{\rm -}\/adic Poincar{\rm \'{\hskip-6pt\it e}}
duality group of dimension {\rm 3,} and if $[K] \in H_3(K,\ZZ_3)$ is a
choice of fundamental class{\rm ,} then
$$
j_\ast[K] \in H_3(S_2^1,\ZZ_3)
$$
is a nonzero $SD_{16}$\/{\rm -}\/invariant generator of infinite order.
\end{lem}

\Proof  The fact that $K$ is a Poincar\'e duality group is
discussed in \cite{Henn}; this discussion is an implementation of the
theory of Lazard \cite{Laz}. We must now address the statements
about $j_\ast[K]$. For this, we first compute with cohomology, and
we use the results and notation of Theorem \ref{thm4.4}.

It is known (see Proposition 4.3 and 4.4 of \cite{Henn}), that the
Lyndon-Serre-Hochschild spectral sequence
$$
H^p(C_3,H^q(K,\FF_3)) \Longrightarrow H^{p+q}(S^1_2,\FF_3)
$$
collapses and that $H^0(C_3,H^q(K,\FF_3))$ is one dimensional for 
$0\leq q\leq 3$; in particular, $j^\ast: H^3(S_2^1,\FF_3) \to H^3(K,\FF_3)$ 
is onto.
Since the composites 
$$
H^1(C_3,\FF_3) \to H^1(S^1_2,\FF_3) \to H^1(w^iC_3w^{-i},\FF_3)
$$
of the inflation with the restriction maps 
are isomorphisms for $i=1,2$, the image of the generator of $H^1(C_3,\FF_3)$
is some linear combination $ax_1 + bx_2$ with both $a \ne 0$ and
$b \ne 0$. This implies that $j^\ast(x_iy_i) = 0$ for $i=1,2$;
for example
$$
aj^\ast(x_1y_1) = j^\ast(y_1(ax_1 + bx_2)) = 0.
$$
But since $j^\ast: H^3(S_2^1,\FF_3) \to H^3(K,\FF_3)$ is onto and 
$H^3(S_2^1,\FF_3)$ is generated by $x_iy_i$ and $a_iy_i$ for $i=1,2$ 
it is impossible that $j^\ast(a_iy_i)$ is trivial for both $i=1,2$. 
Because $K$ is a Poincar\'e duality group of dimension $3$ 
we also know that the Bockstein $\beta:H^2(K,\FF_3) \to H^3(K,\FF_3)$
is zero; hence
$$
j^\ast(a_1y_1 - a_2y_2) = j^\ast(\beta(x_1a_1 - x_2a_2)) = 0 
$$
and therefore 
$$
j^\ast(a_1y_1) = j^\ast(a_2y_2) \ne 0\ .
$$
This shows that $H^3(S_2^1,\FF_3) \to H^3(K,\FF_3)$ is onto and 
factors through the $SD_{16}$-coinvariants, or dually 
$H_3(K,\FF_3) \to H_3(S^1_2,\FF_3)$ is
an injection and lands in the \pagebreak $SD_{16}$-invariants. 
Furthermore, $H_3(K,\FF_3)$ even maps to the kernel of the 
Bockstein $\beta:H_3(S^1_2,\FF_3)\to H_2(S^1_2,\FF_3)$ and the 
induced map 
$$
H_3(K,\FF_3) \to {\frac {\Kern\ \beta: H_3(S^1_2,\FF_3)\to H_2(S^1_2,\FF_3)}
{\Im \beta: H_4(S^1_2,\FF_3)\to H_3(S^1_2,\FF_3)}}\cong 
H_3(S^1_2;\ZZ_3)\otimes_{\ZZ_3}\FF_3
$$ 
is an isomorphism, yielding that $j_\ast[K]$ is a generator of
infinite order. 
\Endproof\vskip4pt
 
We will use cap products with the elements $[K]$ and $j_\ast[K]$ to
construct a commutative diagram for detecting maps $N_3 \to
\ZZ_3[[\GG^1_2/G_{24}]]$. In the form we use the cap product, it
has a particularly simple expression. Let $G$ be a profinite group
and $M$ a continuous $\ZZ_p[[G]]$-module. If $a \in H^n(G,M)$ and
$x \in H_n(G,\ZZ_p)$ we may define $a \cap x \in H_0(G,M)$ as follows:
choose a projective resolution $Q_\bullet \to \ZZ_p$ and represent
$a$ and $x$ by a cocycle $\alpha:Q_n \to M$ and a cycle $y \in \ZZ_p 
\otimes_{\ZZ_p[[G]]} Q_n$ respectively. Then $\alpha$ descends to a
map
$$
\overline{\alpha}: \ZZ_p \otimes_{\ZZ_p[[G]]} Q_n \longrightarrow \ZZ_p 
\otimes_{\ZZ_p[[G]]} M
$$
and $a \cap x$ is represented by $\overline{\alpha}(y)$. 
It is a simple matter
to check that this is well-defined; in particular, if $y = \partial z$
is a boundary, then $\alpha(y) = 0$ because $\alpha$ is a cocycle.
The usual naturality statements apply, which we record in a lemma. 
Note that part (2) is a special case of part (1) (with $K=G$).

\begin{lem}\label{lem4.11} {\rm 1)} If $\varphi:K \to G$ is a 
continuous homomorphism of profinite groups{\rm ,} and 
$x \in H_n(K,\ZZ_p)$ and $a \in H^n(G,M)${\rm ,} then
$$
\varphi_\ast(\varphi^\ast a \cap x) = a \cap \varphi_\ast(x).
$$

{\rm 2)} Suppose $K \subseteq G$ is the inclusion of a normal subgroup and
$M$ is a $G$-module. Then $G/K$ acts on $H^\ast(K,M)$ and $H_\ast(K,\ZZ_p)$
and for $g \in G/K${\rm ,} $a \in H^n(K,M)${\rm ,} and $x \in H^\ast(K,\ZZ_p)$
$$
g_\ast(g^\ast a \cap x) = a \cap g_\ast(x).
$$
\end{lem}

Here is our main diagram. Let $i:K \to \GG^1_2$ be the inclusion.
\vglue-11pt
\begin{equation}\label{eq4.5} {}
\end{equation}
\vglue-46pt
\begin{small}
$$
\xymatrix{
\Hom_{\ZZ_3[[\GG^1_2]]}(N_3,\ZZ_3[[\GG^1_2/G_{24}]]) \rto^{\mathrm{edge}} 
\dto & H^3(\GG^1_2,\ZZ_3[[\GG^1_2/G_{24}]])\dto^{\cap i_\ast[K]} \\
\Hom_{\ZZ_3[SD_{16}]}(\ZZ_3 \otimes_{\ZZ_3[[S^1_2]]} N_3,
\ZZ_3 \otimes_{\ZZ_3[[S^1_2]]} \ZZ_3[[\GG^1_2/G_{24}]])
\rto^-{\mathrm{ev}} \dto &
H_0(\GG_2^1,\ZZ_3[[\GG^1_2/G_{24}]]) \dto \\
\Hom_{\ZZ_3[SD_{16}]}(\FF_3 \otimes_{\ZZ_3[[S^1_2]]} N_3,
\FF_3 \otimes_{\ZZ_3[[S^1_2]]} \FF_3[[\GG^1_2/G_{24}]]) \rto^-{\mathrm{ev}} &
H_0(\GG_2^1,\FF_3[[\GG^1_2/G_{24}]]).\\
}
$$
\end{small}

 \pagebreak\noindent
We now annotate this diagram. The maps labelled $\mathrm{ev}$ are defined
by evaluating a homomorphism at the image of $j_\ast [K]$
under the edge homomorphism
$$
H_3(S_2^1,\ZZ_3) \to \ZZ_3 \otimes_{\ZZ_3[[S^1_2]]} N_3 
$$
resp. 
$$
H_3(S_2^1,\ZZ_3) \to \ZZ_3 \otimes_{\ZZ_3[[S^1_2]]} N_3 \to 
\FF_3 \otimes_{\ZZ_3[[S^1_2]]} N_3
$$
of (\ref{eq4.4}). We now take the image of that element under the projection map to the coinvariants.
Similar remarks apply to $\FF_3$-coefficients.

The diagram commutes, by the definition of cap product.
Theorem \ref{thm4.9} now follows from 
the final two Lemmas \ref{lem4.12} and \ref{lem4.13} below; 
in fact, once we have proved these lemmas, diagram (\ref{eq4.5})
will then show that we can choose a morphism of continuous
$\GG^1_2$-modules
$$
f:N_3 \longrightarrow \ZZ_3[[\GG^1_2/G_{24}]] 
$$
so that
$$
\FF_3 \otimes f: \FF_3 \otimes_{\ZZ_3[[S^1_2]]} N_3 \longrightarrow 
\FF_3 \otimes_{\ZZ_3[[S^1_2]]} \ZZ_3[[\GG^1_2/G_{24}]] 
$$
is nonzero. Then Lemmas \ref{lem4.7}, \ref{lem4.8} and \ref{lem4.3} imply
that $f$ is an isomorphism.
\smallskip


\begin{lem}\label{lem4.12} The homomorphism
$$
\boldsymbol{\cap} i_\ast[K]:H^3(\GG^1_2,\ZZ_3[[\GG^1_2/G_{24}]]) \to
H_0(\GG^1_2,\ZZ_3[[\GG^1_2/G_{24}]])
$$
is an isomorphism.
\end{lem}

\Proof  Recall that we have denoted the inclusion $K\to S_2^1$ 
by $j$. We begin by demonstrating that
$$
\boldsymbol{\cap} j_\ast[K]:H^3(S^1_2,\ZZ_3[[\GG^1_2/G_{24}]]) \to
H_0(S^1_2,\ZZ_3[[\GG^1_2/G_{24}]])
$$
is an isomorphism. Since the action of $C_3$ on $H_3(K,\ZZ_3)\cong \ZZ_3$ 
is necessarily trivial we see that $[K]$ is $C_3$-invariant
and Lemma \ref{lem4.11} supplies a commutative diagram
$$
\xymatrix{
H^3(S^1_2,\ZZ_3[[\GG^1_2/G_{24}]]) \rto^{j^\ast} \dto_{\boldsymbol{\cap} j_\ast[K]} &
H^3(K,\ZZ_3[[\GG^1_2/G_{24}]])^{C_3} \dto^{\boldsymbol{\cap} [K]}\\
H_0(S^1_2,\ZZ_3[[\GG^1_2/G_{24}]]) & \ar[l]_{j_\ast} 
H_0(K,\ZZ_3[[\GG^1_2/G_{24}]])^{C_3}.\\
}
$$
The morphism $\boldsymbol{\cap} [K]$ is an isomorphism by Poincar\'e duality.
As $\ZZ_3[[S^1_2]]$-modules, we have
$$
\ZZ_3[[\GG^1_2/G_{24}]] \cong \ZZ_3[[S^1_2/C_3]] 
\oplus \ZZ_3[[S^1_2/\omega C_3\omega^{-1}]]
$$
on generators $eG_{24}$ and $\omega G_{24}$ in $\GG^1_2/G_{24}$;
hence, as $K$-modules
$$
\ZZ_3[[\GG^1_2/G_{24}]] \cong \ZZ_3[[K]] \oplus \ZZ_3[[K]] 
$$
which shows that $j$ induces an isomorphism on 
$H_0(-,\ZZ_3[[\GG_2^1/G_{24}]])$. 
We claim that $C_3$ acts trivially on $H_0$ and thus 
$j_\ast$ is an isomorphism. In fact, it is clear that the $C_3$-action  
fixes the coset $eG_{24}$; furthermore 
$\omega C_3\omega^{-1}$ is another complement to $K$ in $S_2^1$ 
and therefore 
$$C_3\omega C_3\subset KC_3\omega C_3=S_2^1\omega C_3=
K\omega C_3\omega^{-1}\omega C_3
\subset K\omega C_3 \  , $$
and hence the class of $\omega C_3$  in $H_0$ is also fixed.  

In addition, since $H^q(K,\ZZ_3[[\GG^1_2/G_{24}]]) = 0$ if
$q \ne 3$, the Lyndon-Serre-Hochschild spectral sequence
shows that $j^\ast$ is an isomorphism. 

To finish the proof, we continue in the same manner. Let $r:S^1_2
\to \GG^1_2$ be the inclusion, so that $i = rj:K \to \GG^1_2$. 
By Lemma \ref{lem4.10}, $j_{\ast}[K]$ is $SD_{16}$-invariant 
and then \ref{lem4.11} supplies once more a commutative diagram
$$
\xymatrix{
H^3(\GG^1_2,\ZZ_3[[\GG^1_2/G_{24}]]) \rto^{r^\ast} \dto_{\boldsymbol{\cap} i_\ast[K]} &
H^3(S^1_2,\ZZ_3[[\GG^1_2/G_{24}]])^{SD_{16}} \dto^{\boldsymbol{\cap} j_\ast[K]}\\
H_0(\GG^1_2,\ZZ_3[[\GG^1_2/G_{24}]]) & \ar[l]_{r_\ast} 
H_0(S^1_2,\ZZ_3[[\GG^1_2/G_{24}]])^{SD_{16}}\ .\\
}
$$
We have just shown that $\cap j_\ast[K]$ is an isomorphism.
The map $r_\ast$ sends invariants to coinvariants and, 
since the order of $SD_{16}$
is prime to $3$, is an isomorphism. Again, because the order of $SD_{16}$ 
is prime to $3$ the spectral sequence of the extension 
$S_2^1\to \GG_2^1\to SD_{16}$ collapses at $E_2$ 
and therefore the map $r^\ast$
is an isomorphism. This completes the proof.
\Endproof\vskip4pt 

\begin{lem}\label{lem4.13} The edge homomorphism
$$
\Hom_{\GG^1_2}(N_3,\ZZ_3[[\GG^1_2/G_{24}]]) \to
H^3(\GG^1_2,\ZZ_3[[\GG^1_2/G_{24}]])
$$
is surjective.
\end{lem}

\def\Mod{{{\ZZ_3[[\GG^1_2/G_{24}]]}}}
\Proof  We examine the spectral sequence of 
(4.8): 
$$
E_1^{p,q} \cong \Ext^p_{\ZZ_3[[\GG^1_2]]}(C_q,\Mod) \Longrightarrow
H^{p+q}(\GG^1_2,\Mod).
$$
We need only show that
$$
\Ext^p_{\ZZ_3[[\GG^1_2]]}(C_q,\Mod) = 0
$$
for $p+q = 3$ and $q < 3$. If $q = 1$ or $2$, then $C_q =
\twisted$. Now $\twisted$ is projective as a $\ZZ_3[[\GG_2^1]]$-module and 
therefore $\Ext^p_{\ZZ_3[[\GG^1_2]]}(C_q,\Mod)$ is trivial for $p>0$. 

If $q = 0$, then $C_0 = \Mod$ and by the Shapiro lemma we get an isomorphism 
$$
\Ext^3_{\GG^1_2}(C_0,\Mod) = H^3(G_{24},\Mod) \cong
H^3(C_3,\Mod)^{Q_8}\ .
$$
The profinite $C_3$-set $\GG_2^1/G_{24}$ is an inverse limit of finite 
$C_3$-sets $X_i$ and thus we get an exact sequence 
$$
0\to \lim_i{ }^1H^2(C_3,\ZZ_3[X_i])\to H^3(C_3,\ZZ_3[[\GG_2^1/C_3]])\to 
\lim_iH^3(C_3,\ZZ_3[X_i])\to 0 \ . 
$$ 
Now each $X_i$ is made of a finite number of $C_3$-orbits. 
The contribution of each orbit to $H^3(C_3,-)$ is trivial and to 
$H^2(C_3,-)$ it is either trivial or $\ZZ/3$. 
Therefore $\lim_i$  is clearly trivial and $\lim_i^1$ 
is trivial because the Mittag-Leffler condition is satisfied. 
\hfill\qed

\section{The tower}

In this section we write down the five stage tower whose homotopy
inverse limit is $L_{K(2)}S^0 = E_2^{h\GG_2}$ and the four stage tower
whose homotopy inverse limit is $E_2^{h\GG^1_2}$. As before we 
will write $E = E_2$ and we recall that we have fixed the prime $3$.

To state our results, we will need a new spectrum. Let $\chi$ be the
representation of the subgroup $SD_{16} \subseteq \GG_2$ that appeared in
(\ref{twisted}) and let $e_{\chi}$ be an idempotent in the 
group ring $\ZZ_3[SD_{16}]$ that picks up $\chi$. The action of $SD_{16}$ on 
$E$ gives us a spectrum $E^{\chi}$ which is the telescope 
associated to this idempotent: $E^{\chi}:=e_{\chi}E$. 

Then we have    isomorphisms of Morava modules 
\begin{align*}
E_\ast E^\chi & \cong \Hom_{\ZZ_3[SD_{16}]}(\chi, E_\ast E))\cong 
\Hom_{\ZZ_3[SD_{16}]}(\chi,\Hom^c(\GG_2,E_\ast))\\
& \cong \Hom_{\ZZ_3[[\GG_2]]}(\largetwisted,\Hom^c(\GG_2,E_\ast)) \cong
\Hom^c_{\ZZ_3}(\largetwisted,E_\ast)\ . 
\end{align*}
We recall that 
$\Hom^c(\GG_2,E_{\ast})$ is a Morava module via the diagonal $\GG_2$-\break action,  
and a $\ZZ_3[SD_{16}]$-module via the translation action on $\GG_2$. 
The group\break $\Hom^c_{\ZZ_3}(-,-)$ 
is the group of all homomorphisms which are continuous with 
respect to the obvious $p$-adic topologies. 

It is clear from (\ref{twisted-plus-trivial}) 
that $E^{\chi}$ is a direct summand of 
$E_2^{hQ_8}$ and a  module spectrum over $E_2^{hSD_{16}}$. In fact, 
it is easy to check that the homotopy of $E^\chi$ is free of rank $1$ 
as a $\pi_*(E^{hSD_{16}})$-module on a generator 
$\omega^2u^4\in \pi_8(E^{hQ_8})\subset \pi_*(E)$:   
this generator detemines a map of module spectra 
from $\Sigma^8E^{hSD_{16}}$ to $E^{\chi}$ which is an equivalence. 
From now on we will use this equivalence to replace $E^{\chi}$ by 
$\Sigma^8E^{hSD_{16}}$.  
We note that $E^\chi$ is periodic with period $16$.

\begin{lem}\label{lem5.1} There is an exact sequence of Morava
modules
\begin{align*}
0 \to E_\ast \to E_\ast E^{hG_{24}} &\to
E_\ast \Sigma^8E^{hSD_{16}} \oplus E_\ast E^{hG_{24}}
 \\
&\to  E_\ast \Sigma^8E^{hSD_{16}}  \oplus E_\ast \Sigma^{40}E^{hSD_{16}} \\
& \to E_\ast \Sigma^{40}E^{hSD_{16}}\oplus E_\ast \Sigma^{48}E^{hG_{24}}  
\to E_\ast \Sigma^{48}E^{hG_{24}} \to 0\ .
\end{align*}
\end{lem}

\Proof  Take the exact sequence of continuous $\GG_2$-modules
of Corollary \ref{cor4.2} and apply 
$\Hom^c_{\ZZ_3}(\bullet,E_\ast)$. Then use the isomorphism
$$E_\ast \Sigma^8E^{hSD_{16}}= \Hom_{\ZZ_p}^c(\largetwisted,E_\ast)$$ above 
and the isomorphisms
$$
E_\ast E^{hF} \cong \Hom^c(\GG_2/F,E_\ast)
$$
supplied by (\ref{fixediso}) to get an exact sequence of
Morava modules
\begin{align*}
0 \to E_\ast \to E_\ast E^{hG_{24}}& \to
E_\ast \Sigma^8 E^{hSD_{16}} \oplus  E_\ast E^{hG_{24}}
\to  E_\ast \Sigma^8 E^{hSD_{16}}  \oplus E_\ast \Sigma^8 E^{hSD_{16}} \\
&\to  E_\ast \Sigma^8 E^{hSD_{16}}  \oplus E_\ast E^{hG_{24}} 
\to E_\ast E^{hG_{24}} \to 0\ .
\end{align*}
Finally, we use that $\Sigma^8E^{hSD_{16}}\simeq \Sigma^{40}E^{hSD_{16}}$ 
because $E^{hSD_{16}}$ is periodic of period $16$ 
and $E_\ast E^{hG_{24}} \cong E_\ast \Sigma^{48} E^{hG_{24}}$ as 
Morava modules (see Remark \ref{periodicity}.2).
\phantom{overthere}\hfill\qed

\begin{rem}\label{rem5.2} In the previous lemma, replacing
$\Sigma^8E^{hSD_{16}}$ by $\Sigma^{40}E^{hSD_{16}}$ is merely aesthetic: 
it emphasizes some sort of duality. 
However, $E^{hG_{24}}$ and $\Sigma^{48} E^{hG_{24}}$
are different spectra, even though $E_\ast E^{hG_{24}}
\cong E_\ast \Sigma^{48} E^{hG_{24}}$. This substitution is essential
to the solution to the Toda bracket problem which arises in Theorem
\ref{thm5.5}. 
\end{rem}

In the same way, one can immediately prove

\begin{lem}\label{lem5.3} There is an exact sequence of Morava
modules
$$
0 \to E_\ast E^{h\GG^1_2}\to E_\ast E^{hG_{24}} \to E_\ast \Sigma^8E^{hSD_{16}}
\to E_\ast \Sigma^{40} E^{hSD_{16}} \to E_\ast \Sigma^{48}E^{hG_{24}} \to 0 .
$$
\end{lem}

\begin{thm}\label{thm5.4} The exact sequence of Morava modules
\begin{align*}
0 \to E_\ast \to E_\ast E^{hG_{24}} &\to
E_\ast \Sigma^8E^{hSD_{16}} \oplus E_\ast E^{hG_{24}}
\to  E_\ast \Sigma^8E^{hSD_{16}}  \oplus E_\ast \Sigma^{40}E^{hSD_{16}} \\
& \to E_\ast \Sigma^{40}E^{hSD_{16}}\oplus E_\ast \Sigma^{48}E^{hG_{24}}  
\to E_\ast \Sigma^{48}E^{hG_{24}} \to 0 
\end{align*}
can be realized in the homotopy category of $K(2)$-local spectra
by a sequence of maps 
\begin{align*}
L_{K(2)}S^0 &\to E^{hG_{24}} \to \Sigma^8E^{hSD_{16}} \vee E^{hG_{24}}
\\
&\to  \Sigma^8E^{hSD_{16}} \vee \Sigma^{40} E^{hSD_{16}} \to
\Sigma^{40} E^{hSD_{16}} \vee \Sigma^{48}E^{hG_{24}}
\to \Sigma^{48}E^{hG_{24}}
\end{align*}
so that the composite of any two successive maps is null homotopic.
\end{thm}

\Proof  The map $L_{K(2)}S^0 \to E^{hG_{24}}$ is the
unit map of the ring spectrum $E^{hG_{24}}$. To produce the other
maps and to show that the successive composites are null homotopic, 
we use the diagram of Proposition \ref{prop1.7}. It is enough to show 
that the $E_n$-Hurewicz homomorphism  
$$
\pi_0 F(X,Y) \to \Hom_{E_\ast E}(E_\ast X,E_\ast Y)
$$
is onto when $X$ and $Y$ belong to the set 
$\{\Sigma^8E^{hSD_{16}},E^{hG_{24}}\}$.
(Notice that the other suspensions cancel out nicely, 
since $E^{hSD_{16}}$ is $16$-periodic.) 
Since $\Sigma^8E^{hSD_{16}}$ is a retract of $E^{hQ_8}$, it is
sufficient to show that 
$$
\pi_0 F(E^{hK_1},E^{hK_2}) \to \Hom_{E_\ast E} (E_\ast E^{hK_1},
E_\ast E^{hK_2})
$$
is onto for $K_1$ and $K_2$ in the set $\{ Q_8, G_{24}\}$. 
Using Proposition \ref{prop1.6} and the short exact $\lim$-$\lim^1$ sequence 
for the homotopy groups of $\holim$ we see that it is enough to note that 
$(E_1)^K=0$ and 
$$
\pi_0E^{hK} \to (E_0)^{K}
$$
is surjective whenever $K \subseteq K_2 \cap xK_1x^{-1}$. 
The first part is trivial and for the second part we can appeal to 
Corollary \ref{iso-in-degree-zero}.

To show that the successive compositions are zero, we proceed similarly,
again using Proposition \ref{prop1.7}, but now we have to 
show that various Hurewicz maps are injective. In this case,
the suspensions do not cancel out, and we must show
$$
\pi_0 F(E^{hK_1},\Sigma^{48k}E^{hK_2}) \to \Hom_{E_\ast E} (E_\ast E^{hK_1},
E_\ast \Sigma^{48k}E^{hK_2})
$$
is injective for $k=0$ and $1$, at least for $K_1$ of the form
$\GG_2$, $G_{24}$, or $Q_8$ and $K_2$ of the form $G_{24}$ or $Q_8$.
Since all the spectra involved in Proposition \ref{prop1.6} 
are $72$-periodic, the $\lim$-$\lim^1$ sequence for the 
homotopy groups of $\holim$ shows once more that 
it is sufficient to note that $E_{24k+1}$ is trivial, that 
$\pi_{24k+1}(E^{hK})$ is finite and 
\begin{equation}\label{24k,k=1}
\pi_{24k}E^{hK} \to (E_{24k})^{K}
\end{equation}
is injective for $k=0$ and $1$ and $K \subseteq K_2 \cap xK_1x^{-1}$. 

Again the first part is trivial while the other parts 
follow from Theorem \ref{thm100.13} 
and Corollary \ref{iso-in-degree-zero}.
Note that for $k=1$ the map in (\ref{24k,k=1}) need not be an isomorphism.
\Endproof\vskip4pt 

The following result will let us build the tower.

\begin{thm}\label{thm5.5} In the sequence of spectra
\begin{align*}
L_{K(2)}S^0& \to E^{hG_{24}} \to \Sigma^8E^{hSD_{16}} \vee E^{hG_{24}}
\\
&\to  \Sigma^8E^{hSD_{16}} \vee \Sigma^{40} E^{hSD_{16}} \to
\Sigma^{40} E^{hSD_{16}} \vee \Sigma^{48}E^{hG_{24}}
\to \Sigma^{48}E^{hG_{24}}
\end{align*}
all the possible Toda brackets are zero modulo their indeterminacy.
\end{thm}

\Proof   There are three possible
three-fold Toda brackets, two possible four-fold Toda
brackets and one possible five-fold Toda bracket. All
but the last lie in zero groups.

Because $\Sigma^8E^{hSD_{16}}$ is a summand of $E^{hQ_8}$, 
the three possible three-fold Toda brackets lie in 
\begin{gather*}
\pi_1 F(E^{h\GG_2},E^{hQ_8} \vee\Sigma^{32}E^{hQ_8}),\quad 
\pi_1 F(E^{hG_{24}},\Sigma^{32}E^{hQ_8}\vee \Sigma^{48}E^{hG_{24}}),
\\[4pt]
\pi_1 F(E^{hQ_8} \vee E^{hG_{24}},\Sigma^{48}E^{hG_{24}})
\end{gather*}
which are all zero by Proposition \ref{prop1.6}, Corollary
\ref{zero-in-pi1} and Corollary \ref{effect-of-the-center}.  
The most interesting calculation is
the middle of these three and the most interesting part of that
calculation is
$$
\pi_1 F(E^{hG_{24}}, \Sigma^{48}E^{hG_{24}}) \cong 
\pi_{25}(\holim_i\prod_{G_{24}\backslash\GG_n/U_i)} E^{hH_{x,i}})\ .
$$
This is zero by Corollary \ref{zero-in-pi1} and 
Corollary \ref{effect-of-the-center} (note that the element  
$-1 = \omega^4 \in G_{24}$ is in the center of $\GG_2$ and it is 
in $H_{x,i}$ for every $x$); however, notice that
{\it without} the suspension by $48$ this group is nonzero.

The two possible four-fold Toda brackets lie in
$$
\pi_2 F(E^{h\GG_2},\Sigma^{32}E^{hQ_8}\vee \Sigma^{48}E^{hG_{24}})
\qquad\hbox{and}\qquad
\pi_2 F(E^{hG_{24}},\Sigma^{48}E^{hG_{24}}).
$$
We claim these are also zero groups. All of the calculations
here present some interest. For example, consider
$$
\pi_2 F(E^{hG_{24}}, \Sigma^{48}E^{hG_{24}}) \cong
\pi_{26}(\holim_i\prod_{G_{24}\backslash\GG_n/U_i)} E^{hH_{x,i}})\ .
$$
where $H_{x,i} = xU_ix^{-1} \cap G_{24}$. Since the element
$-1 = \omega^4 \in G_{24}$ is in the center of $\GG_2$,
it is in $H_{x,i}$ for every $x$ and the result follows from 
Corollary \ref{effect-of-the-center} and the observation that 
$\pi_{27}(E^{hH_{x,i}})$ is finite.  

Finally, the five-fold Toda bracket lies in
$$
\pi_3 F(E^{h\GG_2},\Sigma^{48}E^{hG_{24}})\cong\pi_{27}E^{hG_{24}}\cong\ZZ/3\ .
$$
Thus, we do not have the zero group; however, we claim that the map
$E^{h\GG_2}\break \to E^{hG_{24}}$ at the beginning of our sequence supplies a
surjective homomorphism
$$
\pi_\ast F(E^{hG_{24}},\Sigma^{48}E^{hG_{24}})
\to \pi_\ast F(E^{h\GG_2},\Sigma^{48}E^{hG_{24}})\ .
$$
This implies that the indeterminancy of the five-fold
Toda bracket is the whole group, completing the proof.

To prove this claim, note the $E^{h\GG_2} \to E^{hG_{24}}$
is the inclusion of the homotopy fixed points by a larger
subgroup into a smaller one. Thus Proposition~\ref{prop1.6}
yields a diagram

$$
\xymatrix{
F(E^{hG_{24}},E^{hG_{24}}) \rto \dto_{\simeq}
& F(E^{\GG_2},E^{hG_{24}}) \dto^{\simeq}\\
E[[\GG_2/G_{24}]]^{hG_{24}} \rto 
& E[[\GG_2/\GG_2]]^{hG_{24}} 
\rto^-{\simeq} & E^{hG_{24}}\\ 
}
$$
and the contribution of the coset 
$eG_{24}$ in $E[[\GG_2/G_{24}]]^{hG_{24}}$ shows that 
the horizontal map is a split surjection of spectra. 
\Endproof\vskip4pt 

The following result is now an immediate consequence of Theorems
\ref{thm5.4} and \ref{thm5.5}:

\begin{thm}\label{thm5.6} There is a tower of fibrations in the
$K(n)$-local category
\begin{footnotesize}
$$
\hskip-8pt\xymatrix@C=2pt 
{
L_{K(2)}S^0\rto &X_3 \rto & X_2 \rto & X_1 \rto& E^{hG_{24}}\\
\Sigma^{44}E^{hG_{24}}\uto&
\Sigma^{45}E^{hG_{24}}\vee \Sigma^{37}E^{hSD_{16}}\uto&
\Sigma^{6}E^{hSD_{16}}\vee \Sigma^{38}E^{hSD_{16}}\uto&
\Sigma^{7}E^{hSD_{16}}\vee \Sigma^{-1}  E^{hG_{24}}\uto\ .\\
}
$$
\end{footnotesize}
\end{thm}

Using Lemma \ref{lem5.3} and the very same program, we may produce
the following result.  The only difference will be that the Toda
brackets will all lie in zero groups.

\begin{thm}\label{thm5.7} There is a tower of fibrations in the
$K(n)$-local category
$$
\xymatrix{
E^{h\GG^1_2} \rto & Y_2 \rto & Y_1 \rto& E^{hG_{24}}\\
\Sigma^{45}E^{hG_{24}}\uto&\Sigma^{38}E^{hSD_{16}}\uto&
\Sigma^{7}E^{hSD_{16}}\uto\ .\\
}
$$
\end{thm}

\enddocument 

%
%


\references {999}

\bibitem[1]{BL} \name{A. Beauville} and \name{Y. Laszlo},  Un lemme de descente,
{\it C. R. Acad.\ Sci.\ Paris S{\hskip.5pt\rm \'{\hskip-5pt\it e}}r.\ I Math\/}.\ {\bf 320}  (1995),
335--340.

\bibitem[2]{Deligne} \name{P. Deligne}, Courbes elliptiques: formulaire d'apr\`es
{J}. {T}ate, in {\it Modular Functions of One Variable\/}, IV, {\it
Proc.\ Internat.\ Summer School\/} (Univ. Antwerp, Antwerp, 1972), 
53--73, {\it Lecture Notes
in Math\/}.\ {\bf 476}, Springer-Verlag, New York,  1975.

\bibitem[3]{Dev} \name{E. S. Devinatz}, Morava's change of rings theorem, in 
{\it The \v Cech Centennial\/} (Boston, MA, 1993), {\it Contemp.\ Math\/},
{\bf  181}, 83--118, A.\ M.\ S.,
Providence, RI, 1995.

\bibitem[4]{DH} \name{E. S. Devinatz} and \name{M. J. Hopkins}, The action of
the {M}orava stabilizer group on the {L}ubin-{T}ate moduli space of lifts,
{\it Amer.\ J. Math.\/} {\bf 117}  (1995), 669--710.

\bibitem[5]{DH1} \bibline, Homotopy
fixed point spectra for closed subgroups of the Morava stabilizer groups,
{\it Topology} {\bf 43}  (2004), 1--47.

\bibitem[6]{EKMM} \name{A. D. Elmendorf, I. Kriz, M. A. Mandell}, and 
and \name{J. P. May}, {\it Rings\/}, {\it Modules\/}, {\it and Algebras in 
Stable Homotopy Theory}
(With an appendix by M.\ Cole), A.\ M.\ S., Providence, RI, 1997.


\bibitem[7]{GHM} \name{P. Goerss, H.-W.\ Henn}, and \name{M.\ Mahowald},
The homotopy of $L_2V(1)$ for the prime $3$, 
in {\it Categorical Decomposition Techniques in Algebraic Topology\/} 
(Isle of Skye 2001), 125--151, {\it Progr.\ in Math\/}.\  Birkh\"auser,
Boston, 2004. 

\bibitem[8]{GS} \name{V. Gorbounov} and \name{P. Symonds}, 
Toward the homotopy groups of the higher real ${K}$-theory
${E}{\rm {O}}\sb 2$, in {\it Homotopy Theory via Algebraic Geometry and
Group Representations\/} (Evanston, IL, 1997), 103--115, A.\ M.\ S.,
Providence, RI, 1998.


\bibitem[9]{Hart} \name{R. Hartshorne},  {\it Residues and Duality\/},
Lecture notes of a seminar on the work of A.\ Grothendieck,
given at Harvard 1963/64 (With an appendix by P.\ Deligne)
{\it Lecture Notes in Math\/}.\ {\bf 20}, Springer-Verlag, New York, 1966.

\bibitem[10]{Henn} \name{H.-W.\ Henn}, Centralizers of elementary abelian
$p$-subgroups and mod-$p$ cohomology of profinite groups, 
{\it Duke Math.\ J.} {\bf 91} (1998), 561--585.


\bibitem[11]{HG} \name{M. J. Hopkins} and \name{B. H. Gross}, The rigid analytic
period mapping, {L}ubin-{T}ate space, and stable homotopy theory,
{\it Bull.\ Amer.\ Math.\ Soc.\/} {\bf 30}  (1994), 76--86.

\bibitem[12]{HMS} \name{M. J. Hopkins, M. Mahowald,} and \name{H. Sadofsky}, 
Constructions of elements in {P}icard groups, in {\it Topology and
Representation Theory\/} (Evanston, IL, 1992), 89--126, A.\ M.\ S.,
Providence, RI, 1994.







\bibitem[13]{chromatic} \name{M. Hovey}, Bousfield localization functors and
{H}opkins' chromatic splitting conjecture, in {\it The \v Cech
Centennial\/}
(Boston, MA, 1993), {Contemp.\ Math.\/} {\bf 181}, 225--250, A.\ M.\ S.,
Providence, RI, 1995.

\bibitem[14]{HSS} \name{M. Hovey, B. Shipley}, and J\name{. Smith}, Symmetric
spectra, {\it J. Amer.\ Math.\ Soc.\/} {\bf 13}  (2000),  149--208.



\bibitem[15]{HS} \name{M. Hovey} and \name{N. P. Strickland}, {\it  Morava
${K}$-theories and localisation}, {\it Mem.\ Amer.\ Math.\ Soc\/}.\ {\bf
139}, A.\ M.\ S., Providence, R.I. (1999).

\bibitem[16]{Laz} \name{M. Lazard}, Groupes analytiques $p$-adiques,
{\it Inst. Hautes \'Etudes Sci. Publ. Math.\/}.\ {\bf 26} (1965), 389--603.

\bibitem[17]{LT} J. Lubin and J. Tate, Formal moduli for
one-parameter formal {L}ie groups, {\it Bull.\ Soc.\ Math.\ France} {\bf 94}  
(1966), 49--59.

\bibitem[18]{Nave} \name{L. S. Nave}, On the non-existence of Smith-Toda
complexes, Ph.D.\ Thesis, University of Washington, 1999.

\bibitem[19]{Pribble} \name{E. Pribble}, Algebraic stacks for stable homotopy
theory and the algebraic chromatic convergence theorem, Ph.D.\ Thesis,
Northwestern University, 2004.

\bibitem[20]{Rav} \name{D. C. Ravenel}, {\it Complex cobordism and stable
homotopy groups of spheres}, Academic Press Inc., Orlando, FL, 1986.

\bibitem[21]{Nil} \bibline, {\it Nilpotence and Periodicity in
Stable Homotopy Theory} (Appendix C by Jeff Smith), {\it Ann.\ of Math.\
Studies\/} {\bf 128}, Princeton Univ.\ 
Press, Princeton, NJ, 1992.

\bibitem[22]{Rezk} \name{C. Rezk}, Notes on the {H}opkins-{M}iller theorem, in
{\it Homotopy Theory via Algebraic\break Geometry and Group Representations\/}
(Evanston, IL, 1997), 313--366, Amer.\ Math.\  Soc., Providence, RI, 1998.

\bibitem[23]{Shimv} \name{K. Shimomura}, The homotopy groups of the
${L}\sb 2$-localized {T}oda-{S}mith spectrum ${V}(1)$ at the prime $3$,
{\it Trans.\ Amer.\ Math.\ Soc.\/} {\bf 349}  (1997), 1821--1850.

\bibitem[24]{Shim} \bibline, The homotopy groups of the
${L}\sb 2$-localized mod $3$ {M}oore spectrum, {\it J. Math.\ Soc.\ Japan},
{\bf 53}  (2000),  65--90.

\bibitem[25]{ShimW} \name{K. Shimomura} and \name{X. Wang}, 
The homotopy groups $\pi_*(L_2S^0)$ at the prime $3$, 
{\it Topology} {\bf 41}  (2002),  1183--1198.

\bibitem[26]{ShimY} \name{K. Shimomura} and \name{A. Yabe}, The homotopy
groups $\pi\sb *({L}\sb 2{S}\sp 0)$, {\it Topology} {\bf 34} (1995), 
261--289.

\bibitem[27]{Strick} \name{N. P. Strickland}, Gross-Hopkins duality,
{\it Topology} {\bf 39}  (2000), 1021--1033.

\bibitem[28]{SW} \name{P. Symonds} and \name{T. Weigel}, Cohomology of
$p$-adic analytic groups, in
{\it New Horizons in Pro-$p$ Groups}, 349--410, Birkh\"auser, Boston,
2000.


\bibitem[29]{Toda} \name{H. Toda}, Extended $p$-th powers of complexes and
applications to homotopy theory, {\it Proc.\ Japan Acad}.\ {\bf 44} 
(1968), 198--203.

\Endrefs
\end{document}